\documentclass[twocolumn,letter,9pt]{IEEEtran}

\usepackage[utf8]{inputenc} % allow utf-8 input
\usepackage[T1]{fontenc}    % use 8-bit T1 fonts
\usepackage{hyperref}       % hyperlinks
\usepackage{url}            % simple URL typesetting
\usepackage{booktabs}       % professional-quality tables
\usepackage{amsfonts}       % blackboard math symbols
\usepackage{nicefrac}       % compact symbols for 1/2, etc.
\usepackage{microtype}      % microtypography
\usepackage{xcolor}         % colors

\title{Learning of Linear Dynamical Systems \\ as a Non-Commutative Polynomial Optimization Problem}

\author{
    \IEEEauthorblockN{Quan Zhou\IEEEauthorrefmark{1}, Jakub Mare\v{c}ek\IEEEauthorrefmark{3} \thanks{
    An extended abstract entitled: ``Proper Learning of Linear Dynamical Systems as a Non-Commutative Polynomial Optimisation Problem'' by the same authors has been accepted for inclusion in the program of the 24th International Symposium on Mathematical Theory of Networks and Systems. Unfortunately, the conference has been canceled.
    } }\\
    \IEEEauthorblockA{\IEEEauthorrefmark{1}Imperial College London,
    q.zhou22@imperial.ac.uk}\\
    \IEEEauthorblockA{\IEEEauthorrefmark{3}Czech Technical University in Prague,
    jakub.marecek@fel.cvut.cz}
}

\usepackage{times}  % DO NOT CHANGE THIS
\usepackage{helvet} % DO NOT CHANGE THIS
\usepackage{courier}  % DO NOT CHANGE THIS
\usepackage{graphicx} % DO NOT CHANGE THIS
\usepackage{caption} % DO NOT CHANGE THIS AND DO NOT ADD ANY OPTIONS TO IT
\usepackage[utf8]{inputenc} % allow utf-8 input
\usepackage[T1]{fontenc}    % use 8-bit T1 fonts
\usepackage{hyperref}       % hyperlinks
\usepackage{booktabs}       % professional-quality tables
\usepackage{amsfonts}       % blackboard math symbols
\usepackage{nicefrac}       % compact symbols for 1/2, etc.
\usepackage{microtype}      % microtypography

\usepackage{amsmath}
\usepackage{amsthm}
\usepackage{amsfonts}
\usepackage{amssymb}
\usepackage{nccmath}
\usepackage{mathrsfs}

\usepackage{algorithm}
\usepackage{algorithmic}

\usepackage{color}

\usepackage[short]{optidef}

\usepackage{tabularx}
\usepackage{multirow}

\newtheorem{thm}{Theorem}%[section]
\newtheorem{cor}[thm]{Corollary}
\newtheorem{lem}[thm]{Lemma}

\newtheorem{prop}[thm]{Proposition}
\newtheorem{defin}[thm]{Definition}
\newtheorem{defins}[thm]{Definitions}
\newtheorem{assume}[thm]{Assumption}

\newcommand{\RR}{\mathbb{R}}

\newcommand{\Brack}[1]{\left( #1 \right)}

\newcommand{\Exp}[1]{ \mathbb{E} #1}

\newcommand*{\citep}[1]{\cite{#1}}

\DeclareMathOperator*{\std}{std}
\DeclareMathOperator*{\nrmse}{nrmse}
\DeclareMathOperator*{\mean}{mean}
\DeclareMathOperator*{\tr}{tr}
\DeclareMathOperator*{\Sym}{Sym}
\DeclareMathOperator*{\cyc}{cyc}
\DeclareMathOperator*{\II}{II_1}

\DeclareMathOperator*{\spa}{sparse}

\begin{document}

\twocolumn[
\begin{@twocolumnfalse}
\maketitle
\end{@twocolumnfalse}
]

\begin{abstract}
There has been much recent progress in time series forecasting and estimation of system matrices of linear dynamical systems (LDS). 
We present an approach to both problems based on an asymptotically convergent hierarchy of convexifications of a certain non-convex operator-valued problem, which is known as non-commutative polynomial optimization (NCPOP). We present promising computational results,
including a comparison with methods implemented in Matlab System Identification Toolbox.
%We present an approach to both problems based on operator-valued polynomial optimization, which in spite of the non-convexity of the problem, guarantees global convergence of numerical solutions to the least squares estimator. We present promising computational results.
%There has been much recent progress in forecasting the next observation of a linear dynamical system (LDS), which is known as the improper learning, as well as in the estimation of its system matrices, which is known as the proper learning of LDS.
%We present an approach to proper learning of LDS, which in spite of the non-convexity of the problem, guarantees global convergence of numerical solutions to a least squares estimator.
%We present promising computational results.
\end{abstract}

%\begin{keywords}%
%  System Identification; 
%  Linear Dynamic System; Non-Commutative Polynomial Optimization Problem%
%\end{keywords}

\section{Introduction}

We consider the identification of 
vector autoregressive processes with hidden components
from time series of observations, which is a key problem in system identification \cite{Ljung1999}.
Its applications range from the identification of parameters in epidemiological models \citep{anderson1992infectious} and reconstruction of reaction pathways in other biomedical applications \citep{chou2009recent},
to identification of models of quantum systems \cite{bondar2023globally,bondar2022recovering}. 
Beyond this, one encounters either partially observable processes or questions of causality \citep{9363924,geiger2015causal} in almost any application domain.
In the ``prediction-error'' approach to forecasting \citep{Ljung1999},
it allows the estimation of subsequent observations in a time series.

% learning the dynamics of partially observable systems

To state the problem formally, let us define a linear dynamic system $(G,F,V,W)$ as in \cite{WestHarrison}. 
\begin{equation}
\begin{split}
 \phi_{t} &= G \phi_{t-1} + \omega_t, \\ 
Y_t &= F' \phi_t + \nu_t,    
\end{split}
\label{equ:LDS}
\end{equation}
where $\phi_t \in \mathbb{R}^{n\times 1}$ is the hidden state, $Y_t \in\mathbb{R}^{m\times 1}$ is the observed output (measurements, observations), 
$G\in \mathbb{R}^{n\times n}$ and $F\in\mathbb{R}^{n\times m}$ are system matrices, and $\{\omega_t,\nu_t \}_{t\in\mathbb{N}}$ are normally distributed process and observation noises with zero mean and covariance of $W$ and $V$ respectively.
The transpose of $F$ is denoted as $F'$.
Learning (or proper learning) refers to identifying the quadruple $(G,F,V,W)$ given the output $\{Y_t\}_{t\in\mathbb{N}}$. 
We assume that the linear dynamical system $(G,F,V,W)$ is observable \citep{van}, i.e., its observability matrix \citep{WestHarrison}
%:
%\begin{equation}
%\left[ \begin{array}{c}
%F' \\
%F'G \\
%. \\
%. \\
%F'G^{n-1}
%\end{array} 
%\right ]
%\end{equation}
has full rank. Note that a minimal representation is necessarily observable and controllable, cf. Theorem 4.1 in \cite{tangirala2014principles}, so the assumption is not too strong.

%Proper learning is desirable for a number of reasons: outside of the discovery of the actual dynamics, proper learning makes it possible to consider shape constraints easily and to compare the learned model with prior information to check its consistency. For example, the measurement chain has guarantees on the overall error in many applications, which makes it possible to develop guarantees on the magnitude of the measurement noise from above. 

There are three complications. 
First, the dimension $n$ of the hidden state $\phi_t$ is not known, in general. Although \cite{TuhinJMLR} have shown that a lower-dimensional model can approximate a higher-dimensional one rather well, in many cases, it is hard to choose $n$ in practice.
Second, the corresponding optimization problem is non-convex, and guarantees of global convergence have been available only for certain special cases.
Finally, the operators-valued optimization problem is non-commutative, and hence much work on general-purpose commutative non-convex optimization is not applicable without making assumptions \cite[cf.]{bondar2023globally} on the dimension of the hidden state.  

Here, we aim to develop a method for proper learning of LDS that could also estimate the dimension of the hidden state and that would do so with guarantees of global convergence to the best possible estimate, given the observations. 
This would promote explainability beyond what forecasting methods without global convergence guarantees allow for.
%Informally speaking, our approach to proper learning makes it easy to claim ``what you see is what you get'' (because of the convergence guarantees) as well as to explain ``the why'' (because the recovery of the dynamics).
In particular, our contributions are:

\begin{itemize}
    \item We cast learning of a linear dynamical system with an unknown dimension of the hidden state as a non-commutative polynomial optimization problem (NCPOP). This also makes it possible to utilize prior information as shape constraints in the NCPOP.
    \item We show how to use Navascules-Pironio-Acin (NPA) hierarchy \citep{Pironio2010} of convexifications of the NCPOP to obtain bounds and guarantees of global convergence. The runtime is independent of the (unknown) dimension of the hidden state.
   % \item We show that it is possible to use the Gelfand--Naimark--Segal (GNS) construction \citep{gelfand1943imbedding,segal1947irreducible,klep2018minimizer} to extract the optimizers of the NCPOP.
    \item In two well-established small examples of \cite{arima_aaai,hazan2017learning,Jakub}, our approach outperforms standard subspace and least squares methods, as implemented in Matlab\texttrademark{ }  System Identification Toolbox\texttrademark. 
\end{itemize}

%Given a sequence $\{Y_t\}_{t\in\mathbb{N}}$ of output, one can either estimate the subsequent observations, without necessarily learning the system matrices, which is known as improper learning \cite{Ljung1999}, or focus on the recovery of the system matrices, and is known as proper learning.

%
\section{Background}

First, we set our work in the context of related work.
Next, we provide a brief overview of 
non-commutative polynomial optimization, pioneered by 
\cite{Pironio2010} and nicely surveyed by 
\cite{burgdorf2016optimization}, which is our key technical tool.
Prior to introducing our own results, 
we introduce some common notation, %in Section \ref{sec:notation}, 
following \cite{WestHarrison}.% and \cite{Jakub}.

%prediction error approach
%While there is a considerable amount of work on improper learning of linear dynamical systems, proper learning has seen less attention.

%\skipInEightPage{

\subsection{Related Work in System Identification and Control}
\label{sec:relwork}

%Recent research within System Identification appears in venues associated with Control Theory, Statistics, and Machine learning. We refer to \cite{Ljung1999} and \cite{tangirala2014principles} for excellent overviews of the long history of research in the field, going back at least to \cite{ASTROM1965}. In this section, we focus on pointers to key more recent publications. 

There is a long history of research within system identification \cite{Ljung1999}. 
In forecasting under LDS assumptions (improper learning of LDS), a considerable progress has been made in the analysis of predictions for the expectation of the next measurement using auto-regressive (AR) processes in Statistics and Machine Learning. In \cite{anava13}, first guarantees were presented for auto-regressive moving-average (ARMA) processes. In \cite{arima_aaai}, these results were extended to a subset of autoregressive integrated moving average (ARIMA) processes.  \cite{Jakub} have shown that up to an arbitrarily small error given in advance, AR($s$) will perform as well as \emph{any} Kalman filter on any bounded sequence. 
This has been extended by \cite{tsiamis2020online} to Kalman filtering with logarithmic regret.

Another stream of work within improper learning focuses on sub-space methods \citep{katayama2006subspace,van} and spectral methods \cite{hazan2017learning,hazan2018spectral}. 
\cite{tsiamis2019sample,tsiamis2019finite} presented the present-best guarantees for traditional sub-space methods.
\cite{sun2020finite} utilize regularizations to improve sample complexity. 
Within spectral methods, \cite{hazan2017learning} and \cite{hazan2018spectral} have considered learning LDS with input, employing certain eigenvalue-decay estimates of Hankel matrices in the analyses of an auto-regressive process in a %lifted space of 
dimension increasing over time.
We stress that none of these approaches to improper learning  are ``prediction-error'': They do \emph{not} estimate the system matrices.
%in the process of estimating the next measurement.

In proper learning of LDS, many state-of-the-art approaches consider the least squares method, despite complications encountered in unstable systems \citep{faradonbeh2018finite}. \cite{simchowitz2018learning} have  
%that 
provided non-trivial guarantees for the ordinary least squares (OLS) estimator 
in the case of stable $G$ and there being no hidden component, i.e., $F'$ being an identity and $Y_t = \phi_t$. 
Surprisingly, they have also shown that more unstable linear systems are easier to estimate than less unstable ones, in some sense. 
\cite{simchowitz2019learning} extended the results to allow for a certain pre-filtering procedure.
\cite{SarkarRakhlin,TuhinJMLR} extended the results to cover stable,
marginally stable, and explosive regimes.
\cite{oymak2019non} provide a finite-horizon analysis of the Ho-Kalman algorithm.
Most recently, \cite{bakshi2023new} provided a detailed analysis of the use of the method of moments in learning linear dynamical systems, which could be seen as a polynomial-time algorithm for learning a LDS from a trajectory of polynomial length up to a polynomial error.
Our work could be seen as a continuation of the work on the least squares method, with guarantees of global convergence.
%In Computer Science, our work could be seen as an approximation scheme of \cite{vazirani2013approximation},
%as it allows for $\epsilon$ error for any $\epsilon > 0$.

%\cite{lee2019robust}
%\cite{matni2019tutorial}

%Prominent applications of such results have been developed in control theory, esp. in relation to linear quadratic (LQ) regulators, which hare closely related to Kalman filtering once system matrices are known. 
%See  \cite{dean2017sample},  \cite{fazel2018global}, \cite{boczar2019performance}, and most recently, \cite{agarwal2019logarithmic} for prominent examples.

\subsection{Non-Commutative Polynomial Optimization}
\label{sec:ncpop}

Our key technical tool is non-commutative polynomial optimization,
first introduced by \cite{Pironio2010}. 
Here, we provide a brief summary of their results, and refer to \cite{burgdorf2016optimization} for a book-length introduction.
NCPOP is an operator-valued optimization problem with a standard form in \eqref{NCPO}:

\begin{mini}
	  {(\mathcal{H},X,\psi)} {\langle\psi,p(X)\psi\rangle}{}{P^*=}
	  \addConstraint{q_i(X)}{\succcurlyeq 0, }{i=1,\ldots,m}
	  \addConstraint{\langle\psi,\psi\rangle}{= 1,}{}\label{NCPO}
\end{mini}
where $X=(X_1,\ldots,X_n)$ is a tuple of bounded operators on a Hilbert space $\mathcal{H}$.
In contrast to traditional scalar-valued, vector-valued, or matrix-valued optimization techniques, the dimension of variables $X$ is unknown \textit{a priori}.
Let $[X,X^{\dag}]$ denotes these $2n$ operators, with the $\dag$-algebra being conjugate transpose.
The normalized vector $\psi$, i.e., $\|\psi \|^2=1$ is also defined on $\mathcal{H}$ with the inner product $\langle\psi,\psi\rangle=1$. 
$p(X)$ and $q_i(X)$ are polynomials 
and $q_i(X)\succcurlyeq 0$ denotes that the operator $q_i(X)$ is positive semidefinite. 
Polynomials $p(X)$ and $q_i(X)$ of degrees $\deg(p)$ and $\deg(q_i)$, respectively, can be written as:
\begin{equation}
    p(X)=\sum_{|\omega|\leq \deg(p)} p_{\omega} \omega,\quad
    q_i(X) = \sum_{|\mu|\leq \deg(q_i)} q_{i,\mu} \mu, 
    \label{LCoP}
\end{equation}
where $i = 1,\ldots,m$. Monomials $\omega,\mu,u$ and $\nu$ in following text are products of powers of variables from $[X,X^{\dag}]$.
The degree of a monomial, denoted by $|\omega|$, refers to the sum of the exponents of all operators in the monomial $\omega$.
Let $\mathcal{W}_k$ denote the collection of all monomials whose degrees $|\omega|\leq k$, or less than infinity if not specified.
Following \cite{akhiezer1962some}, we can define the moments on field $\mathbb{R}$ or $\mathbb{C}$, with a feasible solution $(\mathcal{H},X,\psi)$ of problem \eqref{NCPO}:
\begin{equation}
    y_{\omega} = \langle \psi, \omega\; \psi \rangle,\label{DEFoMOMENT}
\end{equation}
for all $\omega\in\mathcal{W}$ and $y_1=\langle \psi,\psi \rangle=1$.
Given a degree $k$, the moments whose degrees are less or equal to $k$ form a sequence $y=(y_{\omega})_{|\omega| \leq 2k}$.
We call $k$ as the moment order.
With a finite set of moments $y$ of moment order $k$, we can define a corresponding $k^{th}$-order moment matrix $M_k(y)$:
\begin{equation}
    M_k(y)(\nu,\omega) = y_{\nu^{\dag}\omega} = \langle \psi, \nu^{\dag}\;\omega\; \psi \rangle, \label{equ:MomentMatrix}
\end{equation}
for any $ |\nu|,|\omega| \leq k$ and the localizing matrix $M_{k_i}(q_i y)$: 
\begin{align}
    M_{k_i}(q_iy)(\nu,\omega) & = \sum_{|\mu| \leq \deg(q_i)} q_{i,\mu} y_{\nu^{\dag}\mu\omega} \label{equ:LocalizingMatrix}\\ & = \sum_{|\mu| \leq \deg(q_i)} q_{i,\mu} \langle \psi, \nu^{\dag}\; \mu\; \omega\; \psi \rangle, \notag
\end{align}
for any $|\nu|,|\omega| \leq k_i$, where $k_i=k-\lceil \deg(q_i)/2\rceil$, and $i=1,\dots,m$. %The upper bounds of $|\nu|$ and $|\omega|$ are lower than that of moment matrix because $y_{\nu^{\dag}\mu\omega}$ is only defined on $\nu^{\dag}\mu\omega \in \mathcal{W}_{2k}$ while $\mu\in \mathcal{W}_{\deg(q_i)}$.

%\textcolor{red}{Prune?}

If $(\mathcal{H},X,\psi)$ is feasible, 
one can utilize the Sums of Squares theorem of \cite{helton2002positive} and \cite{mccullough2001factorization} to derive semidefinite programming (SDP) relaxations.
In particular, we can obtain a $k^{th}$-order SDP relaxation of the non-commutative polynomial optimization problem \eqref{NCPO} by choosing a moment order $k$ that satisfies the condition of $2k\geq \max\{\deg(p), \deg(q_i)\}$.
%Recall that $\deg(p)$ denotes the degree of $p(X)$ and $\deg(q_i)$ denote the degree of $q_i(X)$ for $i=1,\dots,m$. 
In the Navascules-Pironio-Acin (NPA) hierarchy \citep{Pironio2010}, the SDP relaxation of moment order $k$, has the following form:
\begin{mini}
{ y=(y_{\omega})_{|\omega|\leq 2k} }{\sum_{|\omega|\leq d} p_{\omega} y_{\omega}}{}{P^k=}
\addConstraint{M_k(y)}{\succcurlyeq 0}{}
\addConstraint{M_{k_i}(q_i y)}{\succcurlyeq 0, }{i=1,\ldots,m}
\addConstraint{y_1}{= 1.}{}
%\addConstraint{\langle\psi,\psi\rangle}{= 1.}{}
\label{NCPO-R}
\end{mini}
Notice that there are variants \cite{wang2019tssos,wang2020chordal,wang2020sparsejsr} that exploit the sparsity and significantly reduce the computational burden. 

Let us define the quadratic module, following \cite{Pironio2010}.  Let $Q=\{q_i,i=1,\dots,m\}$ be the set of polynomials determining the constraints. 
%defining $\Upsilon$ (\ref{polyncprog2}). 
The \emph{positivity domain} $\mathbf{S}_Q$ of $Q$ are $n$-tuples of bounded operators  $X=(X_1,\ldots,X_n)$ on a Hilbert space $\mathcal{H}$ making all $q_i(X)$ positive semidefinite.
The \emph{quadratic module} $\mathbf{M}_Q$ is the set of  $\sum_if_i^{\dag}f_i+\sum_i\sum_j g_{ij}^{\dag}q_ig_{ij}$ 
where $f_i$ and $g_{ij}$ are polynomials from the same ring. %in $\mathbb{K}[x,x^*]$.
As in \cite{Pironio2010}, we assume:

\begin{assume}[Archimedean]
\label{Archimedean}
Quadratic module $\mathbf{M}_Q$ of \eqref{NCPO} is Archimedean, i.e., there exists a real constant $C$ such that $C^2-(X_1^{\dag}X_1+\cdots+X_{2n}^{\dag}X_{2n})\in \mathbf{M}_Q$. 
\end{assume}

%We refer to the Supplementary material for further discussion.
If the Archimedean assumption is satisfied, Pironio et al.\ \cite{Pironio2010} have shown that $\lim_{k \to \infty} P^k=P^*$
and how to use the so-called rank-loop condition \cite{Pironio2010} to detect global optimality.
We refer to an extended version online \cite{zhou2020proper} for further details.

\subsection{Minimizer Extraction and Gelfand-Naimark-Segal Construction}
\label{sec:GNS}

Notice that the solution of the SDP relaxation makes it possible to read out the value of the objective function $\langle\psi,p(X)\psi\rangle$ of \eqref{NCPO} 
easily, by looking up the correct entries of the moment matrix \eqref{equ:MomentMatrix}.
To extract the optimizer with this objective-function value, one may utilize a variant of the singular-value decomposition of the moment matrix 
pioneered by \cite{henrion2005detecting},
which can be construed \cite{klep2018minimizer} as the Gelfand--Naimark--Segal (GNS) construction \citep{gelfand1943imbedding,segal1947irreducible,dixmier1969algebres}.
(The GNS construction essentially produces a *-representation from a positive linear functional of a C*-algebra on a Hilbert space. Under the Archimedean assumption, this method could be applied to non-commutative polynomials, which are not C*-algebras otherwise.
We refer to an extended version online \cite{zhou2020proper} for further details.)
These SVD-based approaches do not require the rank-loop condition to be satisfied, 
as is well explained in Section 2.2 of \cite{klep2018minimizer}.
Once global optimality is detected (cf. the previous section), it is possible to extract the global optimum $(H^*, X^*, \psi^*)$ from the solution of the SDP relaxation of \eqref{NCPO} by Gram decomposition; cf. Theorem 2 in \cite{Pironio2010}.

%\subsection{Notation and Definitions}
%\label{sec:notation}

%The central quantity within system identification is the estimate of the next observation, given the current data, i.e., the expectation of $Y_t$ conditional on $Y_{1} \dots Y_{t-1}$: 
%\begin{equation}
%f_{t+1} := \Exp{\Brack{Y_{t+1}\mid  Y_{1} \dots Y_{t-1}}} := F'a_{t+1}.\label{f_t+1}
%\end{equation}
%Traditionally, one writes down the recursive formula:
%\begin{eqnarray}
%a_{t+1} &=& Gm_{t} = GA_t Y_t + Z_t a_t \nonumber\\ 
%&=& GA_t Y_t + Z_t G A_{t-1} Y_{t-1} + Z_t Z_{t-1} a_{t-1} \label{eq:a_t_update} \\ 
%&\dots \nonumber
%\end{eqnarray}
%from the recursive update equations of Kalman filter (cf. the Supplementary material),
%where $Z_t$ denotes $G(I-F\otimes A_t)$. 
%By unrolling the recursive update equations of Kalman filter, $a_{t+1}$ can be expressed in terms of $a_{t-s}$ and observations %$Y_{t-s},\dots,Y_{t}$. Hence, substituting $a_{t+1}$ in \eqref{f_t+1} by the recursive formula, %\eqref{eq:a_t_update}, 
%\cite{Jakub} extend the forecast of $Y_{t+1}$ given $Y_1,\dots,Y_t$ to an auto-regressive model with the degree of $s+1$ plus a remainder term.
%As we detail in the Supplementary material, 
%This approach also yields a non-commutative polynomial optimization problem, but of a considerable degree.
%Our main result hence considers a rather different approach.

%as in Figure \ref{fig:exp}.

%Since the remainder term goes to zero exponentially fast with $s$ \cite{Jakub}, we can ignore the remainder term when using Figure \ref{fig:exp} for forecast.

\section{The Main Result}
\label{sec:main}

Given a trajectory of observations $Y_1$,...,$Y_{t-1}$, 
loss is a one-step error function at time $t$ that compares an estimate $f_t$ with the actual observation $Y_t$. 
Within the least squares estimator, we aim to minimize the sum of quadratic loss functions, i.e., 
\begin{equation}
\label{obj}
    \min_{f_{t,t\geq 1}} \sum_{t\geq 1} \|Y_{t}-f_{t}\|^2,
\end{equation}
where the estimates $f_t,t\geq 1$ are decision variables.
The properties of the optimal least squares estimate are well understood: it is consistent, cf. Mann and Wald \cite{mann1943statistical} and Ljung \cite{ljung1976consistency}, and has favorable sample complexity, cf. Theorem 4.2 of Campi and Weyer \cite{campi2002finite} in the general case, and to Jedra and Proutiere \cite{jedra2020finite} for the latest result parameterized by the size of a certain epsilon net.
%(capturing the intrinsic complexity, similar to VC dimension).
We stress, however, that \emph{it has not been understood} how to solve the non-convex optimization problem, in general, outside of some special cases \cite{Hardt} and recent, concurrent work of \cite{bakshi2023new}. In contrast to \cite{Hardt}, we focus on a method achieving global convergence under mild assumptions, and specifically without assuming the dimension of the hidden state is known. 

When the dimension of the hidden state is not known, we need operator-valued variables $m_t$ to model the state evolution, and some additional scalar-valued variables. We denote the process noise and the observation noise at time $t$ by $\omega_t$ and $\nu_t$, respectively. We also denote as such  
the decision variables corresponding to the estimates thereof, if there is no risk of confusion.   
%Because $Y_t$ is already defined to be a scalar, observation noise $\nu_t$ is also a scalar. 
%Since the process and observation noises are normally distributed with mean zero (by assumption), we would expect a sum of $\omega_t$ and a sum of $\nu_t$ over $t$ to be close to $0$, for sufficiently large number of summands. Therefore, 
If we add the sum of the squares of $\omega_t$ and the sum of the squares of $\nu_t$ as regularizers to the objective function with sufficiently large multipliers and minimize the resulting objective, we %can try different values of $B$ 
should reach a feasible solution with respect to the system matrices with the process noise $\omega_t$ and observation noise $\nu_t$ being close to zero.

%In this model, the observations are from the linear system which does not have process-noise and observation-noise, so we expect $p_t$ to be as close to $0$ as possible. If we add the sum of the squares of $p_t$ to the objective with a multiplier $B$ and minimize the sum, we can try different values of $B$ to reach a feasible solution with the noise $p_t$ being close to zero.

Overall, such a formulation has the  form in Equations~\eqref{obj_B} subject to (\ref{NFF_1}--\ref{NFF_2}).
The inputs are $Y_{t},t\geq 1$, i.e., the time series of the actual measurements, of a time window $T$ thereof, and multipliers $c_1, c_2$.
Decision variables are system matrices $G$, $F$; noisy estimates $f_t$, realizations $\omega_t$, $\nu_t$ of noise, for $t\geq 1$;
and state estimates $m_t$, for $t\geq 0$, which include the initial state $m_0$. 
We minimize the objective function:
\begin{equation}
\min_{f_t, m_t, G, F, \omega_t, \nu_t} \sum_{t\geq 1} \| Y_{t}-f_{t} \|^2  + c_1 \sum_{t\geq 1} \nu_t^2 + c_2 \sum_{t\geq 1} \omega_t^2 \label{obj_B} 
\end{equation}
for a 2-norm $\| \cdot \|$ over the feasible set given by constraints for $t\geq 1$:
\begin{align}
m_t & = G m_{t-1} + \omega_t \label{NFF_1} \\
f_{t} & = F' m_t + \nu_t. \label{NFF_2}
\end{align}
We call the term $F'm_t$ noise-free estimates, which are regarded as our simulated/ predicted outputs.
Equations~\eqref{obj_B} subject to (\ref{NFF_1}--\ref{NFF_2}) give us the least squares model. We can now apply the techniques of non-commutative polynomial optimization to the model so as to recover the system matrices of the underlying linear system.

\begin{thm}
\label{T1}
For any observable linear system $(G,F,V,W)$, 
for any length $T$ of a time window,
and any error $\epsilon > 0$, 
under Assumption \ref{Archimedean},
there is a convex optimization problem whose objective function value is at most $\epsilon$ away from \eqref{obj_B} subject to (\ref{NFF_1}--\ref{NFF_2}). 
Furthermore, an estimate of $(G,F,V,W)$ can be extracted from the solution of the same convex optimization problem.
%  the best possible estimate of system matrices of the system $L$ based on the $T$ observations
%, via the least squares model, 
\end{thm}

%To do so, we need to have the non-commutative versions of those formulations firstly, thus introduce a a Hilbert space $\mathcal{H}$ with all operators and a normalized vector defined on $\mathcal{H}$. 

%Proof is in the Supplementary material.

\begin{proof}
First, we need to show the existence of a sequence 
of convex optimization problems, whose objective function approaches the optimum of the non-commutative polynomial optimization problem.
As explained in Section \ref{sec:ncpop} above, 
\cite{Pironio2010} show that there is a sequence of natural semidefinite-programming relaxations of \eqref{NCPO}.
The convergence of the sequence of their objective-function values is shown by Theorem 1 of \cite{Pironio2010}, which requires Assumption \ref{Archimedean}.
The translation of a problem involving multiple scalar- and operator-valued variables $f_t, m_t, G, F, \omega_t, \nu_t$ in (\ref{obj_B}--\ref{NFF_2}) to $(\mathcal{H},X,\psi)$ of \eqref{NCPO}, also known as the product-of-cones construction, is somewhat tedious, but routine and implemented in multiple software packages \cite[e.g.]{wittek2015algorithm,wang2021exploiting}.
Second, we need to show that extraction of an estimate of $(G,F,V,W)$ from the SDP relaxation of order $k(\epsilon)$ in the series is possible. 
There, one utilizes the Gelfand--Naimark--Segal (GNS) construction \citep{gelfand1943imbedding,segal1947irreducible}, as explained in Section 2.2 of \cite{klep2018minimizer}
or in Section \ref{sec:GNS} above. Notice that \cite[cf.]{lee2023computability} the estimate of $(G,F,V,W)$ may have a higher error than $\epsilon$.  
\end{proof}

%\subsection{Extensions}

%Notice that there are several other approaches possible, as discussed in the Supplementary material, and the same reasoning can be applied there as well.
This reasoning can be applied to more 
complicated formulations, involving shape constraints.
%For example, consider symmetric system matrices. 
%A complex square matrix $X$ is called Hermitian if $X^{\dag} = X$ and two variables $X_1$ and $X_2$ commute if $[X_1,X_2]=X_1 X_2 - X_2 X_1=0$, where the operator $[\cdot,\cdot]$ is the commutator of two elements.
For instance, in quantum systems \cite{bondar2023globally}, density operators are Hermitian and this constraint can be added to the least squares formulation.

Crucially for the practical applicability of the method, one should like to exploit the sparsity in the NCPOP (\ref{obj_B}--\ref{NFF_2}). Notice that one can decompose the problem (\ref{obj_B}--\ref{NFF_2}) into $t$ subsets of variables involving 
$f_t, m_{t-1}, m_t, G, F, \omega_t, \nu_t$, which satisfy the so-called running intersection property \cite{wang2021exploiting}. We refer to \cite{klep2019sparse} for a seminal paper on trace optimization exploiting correlative sparsity, and to \cite{wang2021exploiting} for the variant exploiting term sparsity. We also present a brief summary online \cite{zhou2020proper}.

Also, note that the extraction of the minimizer using the GNS procedure, as explained in Section \ref{sec:GNS} above, is stable to errors in the moment matrix, for \emph{any} NCPOP, including the pre-processing above.
See Theorem 4.1 in \cite{klep2018minimizer}.
That is: it suffices to solve the SDP relaxation with a fixed error, in order to extract the minimizer. 

One can also utilize a wide array of reduction techniques on the resulting SDP relaxations. Notable examples include facial reduction \citep{borwein1981facial,permenter2018partial} and exploiting sparsity \citep{fukuda2001exploiting}. Clearly, these can be applied to any SDPs, irrespective of the non-commutative nature of the original problem, but can also introduce \citep{kungurtsev2018two} numerical issues.
We refer to \cite{majumdar2019recent} for an up-to-date discussion.

\section{Numerical illustrations}
\label{app:numerical}

Let us now present the implementation of the approach using the techniques of non-commutative polynomial optimization  \citep{Pironio2010,burgdorf2016optimization} and to compare the results with traditional system identification methods.
Our implementation is available online \footnote{\url{https://github.com/Quan-Zhou/Proper-Learning-of-LDS}}.
%In our implementation, we make use of a globally convergent Navascués-Pironio-Acín (NPA) hierarchy \citep{Pironio2010} of semidefinite programming (SDP) relaxations, as utilized in the proof of Theorem~\ref{T1}, and its sparsity-exploiting variant, known as the term-sparsity exploiting moment/SOS (TSSOS) hierarchy \citep{wang2019tssos,wang2020chordal}.
%Because the degrees of objective \eqref{obj_B} and constraints in (\ref{NFF_1}--\ref{NFF_2}) are all less or equal to $2$, the moment order $k$ within the respective hierarchy can start from $k=1$. %and increase by $1$ in each iteration. 
%The SDP of a given order in the respective hierarchy is constructed using \texttt{ncpol2sdpa 1.12.2}\footnote{\url{https://github.com/peterwittek/ncpol2sdpa}} of Wittek \citep{wittek2015algorithm} or the \texttt{tssos}\footnote{\url{https://github.com/wangjie212/TSSOS}} of Wang et al. \citep{wang2019tssos,wang2020chordal} 
%and solved by \texttt{mosek 9.2} or \texttt{sdpa} of Yamashita et al. \citep{yamashita2003implementation}.\textcolor{red}{how to cite mosek} 
%Empirically, we plot the run-time as a function of $T$ in Figure~\ref{fig:runtime}.
%Our implementation is available online for review purposes and will be open-sourced upon acceptance. 
%It relies on \cite{wittek2015algorithm}, \cite{cafuta2011ncsostools}, and SDPA of \cite{yamashita2003implementation}. As comparison, we also present the runtime of SDP of exploiting term-sparsity (TSSOS), as pioneered by Wang et al. \cite{wang2019tssos,wang2020chordal}.
We present our experimental settings in more detail online \cite{zhou2020proper}.

\subsection{The general setting}
\label{sec:settings}

%Consider a 2-dimensional system ($n=2$), with$G=\bigl(\begin{smallmatrix} 0.99&0\\1&0.2 \end{smallmatrix}\bigr)$, $F'=\bigl(\begin{smallmatrix} 1&0.8\end{smallmatrix}\bigr)$ and the starting point $m_0'=\bigl(\begin{smallmatrix} 1 & 1 \end{smallmatrix}\bigr)$. 

\paragraph{Our formulation and solvers}

For our formulation, we use Equations~\eqref{obj_B} subject to (\ref{NFF_1}--\ref{NFF_2}), where we need to specify the values of $c_1$ and $c_2$.
To generate the SDP relaxation of this formulation as in \eqref{NCPO-R}, we need to specify the moment order $k$. Because the degrees of objective \eqref{obj_B} and constraints in (\ref{NFF_1}--\ref{NFF_2}) are all less than or equal to $2$, the moment order $k$ within the respective hierarchy can start from $k=1$.

In our implementation, we use a globally convergent Navascués-Pironio-Acín (NPA) hierarchy \citep{Pironio2010} of SDP relaxations, as utilized in the proof of Theorem~\ref{T1}, and its sparsity-exploiting variant, known as the non-commutative variant of the term-sparsity exploiting moment/SOS (TSSOS) hierarchy \citep{wang2019tssos,wang2020chordal,wang2021exploiting}.
(See \cite{zhou2020proper} for a summary.)
The SDP of a given moment order within the NPA hierarchy is constructed using \texttt{ncpol2sdpa} 1.12.2\footnote{\url{https://github.com/peterwittek/ncpol2sdpa}} of Wittek \citep{wittek2015algorithm}.
The SDP of a given moment order within the non-commutative variant of the TSSOS hirarchy is constructed using the \texttt{nctssos}\footnote{\url{https://github.com/wangjie212/NCTSSOS}} of Wang et al. \citep{wang2021exploiting}.
%with trace minimization implemented.
Both SDP relaxations are then solved by \texttt{mosek} 9.2 \cite{mosek2020mosek}.
%Note that our method is implemented using NPA hierarchy unless specified.

\paragraph{Baselines}
We compare our method against leading methods for estimating state-space models, as implemented in MathWorks\texttrademark{ } Matlab\texttrademark{ } System Identification Toolbox\texttrademark. Specifically, we test against a combination of least squares algorithms implemented in routine \texttt{ssest}  (``least squares auto''), subspace methods of \cite{van} implemented in routine \texttt{n4sid} (``subspace auto''), and a subspace identification method of \cite{jansson2003subspace} with an ARX-based algorithm to compute the weighting, again utilized via \texttt{n4sid} (``ssarx''). 

To parameterize the three baselines, we need to specify the dimension $d$ of the estimated state-space model. We would set $d=n$ directly or alternatively, iterate from $1$ to the highest number allowed in the toolbox when the underlying system is unknown, e.g., in real-world stock-market data.
Then, we need to specify the error to be minimized in the loss function during estimation.
In fairness to the baselines, we use the one-step ahead prediction error when comparing prediction performance and simulation error between measured and simulated outputs when comparing simulation performance.

\paragraph{The performance index}
To measure the goodness of fit between the ground truth $\{Y_t\}_{t=1}^{T}$ (actual measurements) and the noise-free simulated/ predicted outputs $\{F'm_t\}_{t=1}^{T}$, using different system identification methods, we introduce the
\textit{normalized root mean squared error (nrmse)} fitness value: 
\begin{equation}
\mathrm{\nrmse}:=\left(1- \frac{\left\|Y-F'm\right\|^2}{\left\|Y-\mean(Y)\right\|^2}\right)\times 100\%,\label{NRMSE}
\end{equation}
where $Y$ and $F'm$ are the vectors consisting of the sequence $\{Y_t\}_{t=1}^{T}$ and $\{F'm_t\}_{t=1}^{T}$ respectively.
A higher nrmse fitness value indicates better simulation or prediction performance. 

\subsection{Experiments on the example of Hazan et al.}
\label{sec:exp:hazan}

Experiments in Sections \ref{sec:exp:hazan}--\ref{sec:exp:baselines} utilize synthetic time series of $T$ observations generated using LDS of the form in \eqref{equ:LDS}, with the tuple $(G,F,V,W)$ and the initial hidden state $\phi_0$ 
 detailed next.
We use the dimension $n$ to indicate that the time series of observations were generated using $n \times n$ system matrices,
while we use operator-valued variables to estimate these. 
The standard deviations of process noise and observation noise $W,V$ are chosen from $0.1,0.2,\dots,0.9$. 
Note that $W$ is an $n\times n$ matrix in general, while we consider the spherical case of $W=0.1\times I_{n}$, where $I_n$ is the $n$-dimensional identity matrix, which we denote by $W=0.1$.
%\textcolor{red}{Let $T$ be the length of time window. }
%We would take a sequence of measured output from the LDS as the time series $\{Y_t\}_{t=1}^{T}$, with $T$ be the length of the time window.

In our first experiment, we explore the statistical performance of feasible solutions of the SDP relaxation using the example of Hazan et al.\ \cite{hazan2017learning,Jakub}.
We performed one experiment on each combination of standard deviations of process $W$ and observation noise $V$ from the discrete set $0.1,0.2,\dots,0.9$, i.e., 81 runs in total.
%For our method, we set parameters $c_1=5\times 10^{-4},c_2=10^{-4}$, and the moment order $k=1$, then the formulation is relaxed via \texttt{ncpol2sdpa} 1.12.2 and solved by \texttt{mosek} 9.2.
%, in comparison with other  system identification methods.
%In order to generate the sequence of observations, we build a linear system (i.e., the ground true) with the underlying dynamic being 2-dimensional. 

\begin{figure}[t!]
\centering
\includegraphics[width=0.35\textwidth]{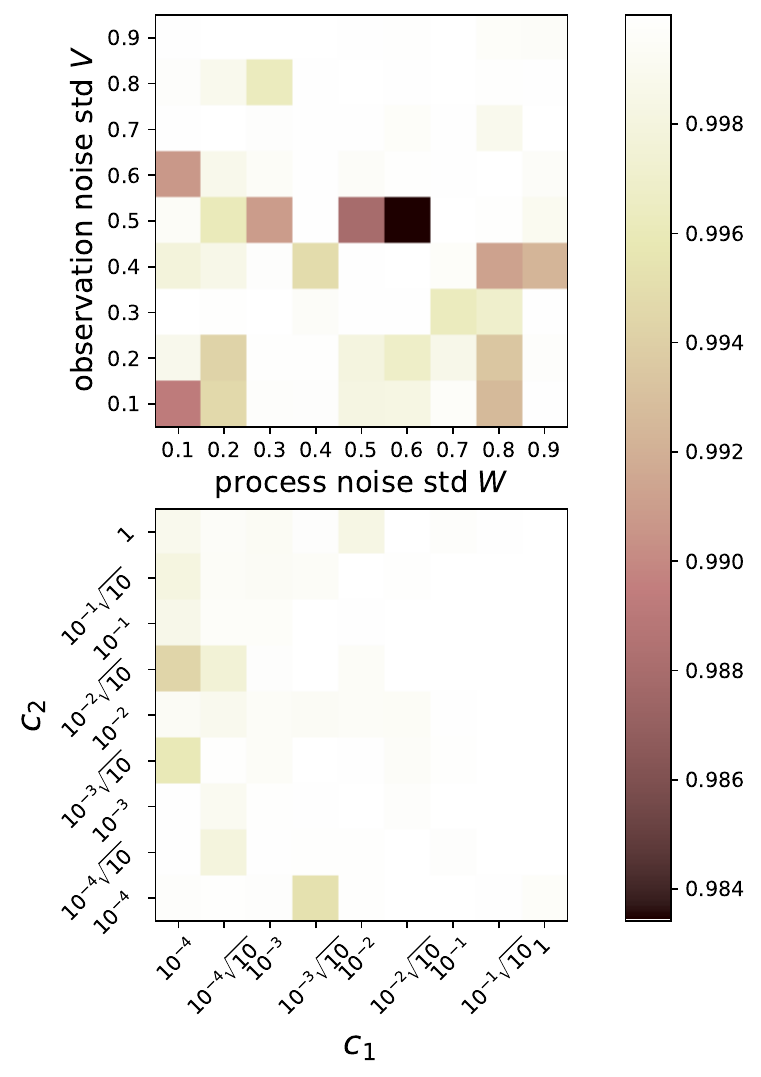}
\caption{\textbf{Upper:} The nrmse fitness values \eqref{NRMSE} of $81$ experiments of our method at different combinations of noise standard deviations of process noise $W$ and observation noise $V$ and \textbf{Lower:} at different combinations of parameters $c_1$ and $c_2$. Both use the data generated from systems in \eqref{equ:LDS}.
Lighter colors indicate higher nrmse and thus better simulation performance.
}\label{fig:NCPO100}
\end{figure}

Figure~\ref{fig:NCPO100} illustrates the nrmse values of the $81$ runs of our method in different combinations of standard deviations of process noise $W$ and observation noise $V$ (upper), and another $81$ experiments in different combinations of $c_1$ and $c_2$ (lower). 
In the upper subplot of Figure~\ref{fig:NCPO100}, we consider: $n=2$, $G=\bigl(\begin{smallmatrix} 0.9&0.2\\0.1&0.1\end{smallmatrix}\bigr)$, $F'=\bigl(\begin{smallmatrix} 1&0.8\end{smallmatrix}\bigr)$, the starting point $\phi_0'=\bigl(\begin{smallmatrix} 1 & 1 \end{smallmatrix}\bigr)$, and $T=20$. 
In the lower subplot of Figure~\ref{fig:NCPO100}, we have the same settings as in the upper one, except for $W=V=0.5$ and the parameters $c_1,c_2$ being chosen from $10^{-4},\dots,1$. %We further perform $9\times 9$ experiments on each combination of $c_1$ and $c_2$.
%are illustrated in the upper subplot of Figure~\ref{fig:NCPO100}.
%Note that we only use the feasible solutions of the relaxation problems, which may explain some of the ``non-linear'' nature of the plot.
It seems clear the highest nrmse is to be observed for the standard deviation of both process and observation noises close to $0.5$.
While this may seem puzzling at first, notice that higher standard deviations of noise make it possible to approximate the observations by an auto-regressive process with low regression depth \cite[Theorem 2]{Jakub}. 
The observed behavior is therefore in line with previous results \cite[e.g., Figure 3]{Jakub}. 

\subsection{Comparisons against the baselines}
\label{sec:exp:baselines}

Next, we investigate the simulation performance of our method in comparison with other system identification methods, for varying LDS used to generate the time series. 
Our method and the three baselines described in Section \ref{sec:settings} are run $30$ times for each choice of the standard deviations of the noise, %i.e., $4\times 30\times 9$ runs in total,
with all methods using the same time series.

\begin{figure}[t!]
\centering
\includegraphics[width=0.35\textwidth]{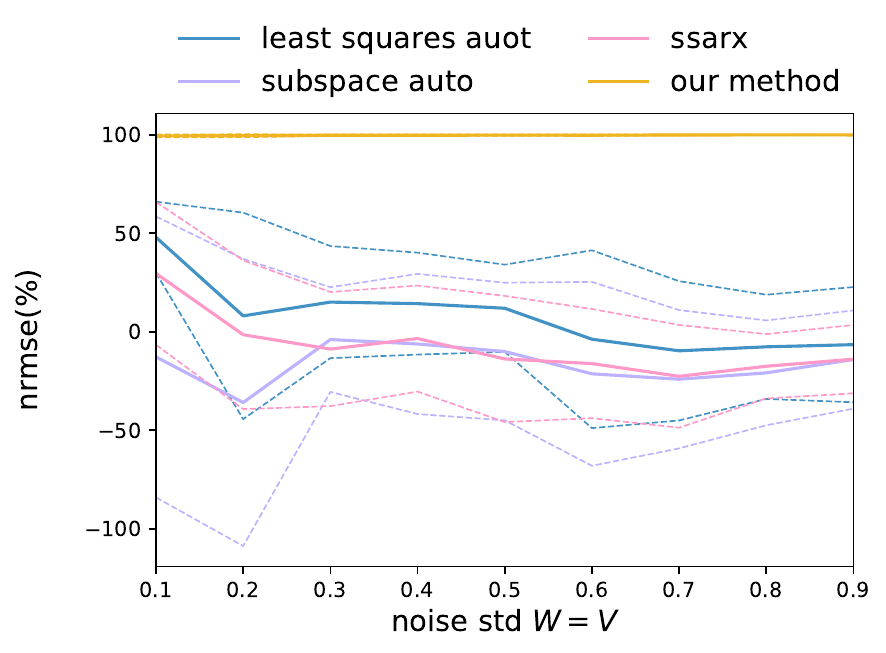}
\includegraphics[width=0.35\textwidth]{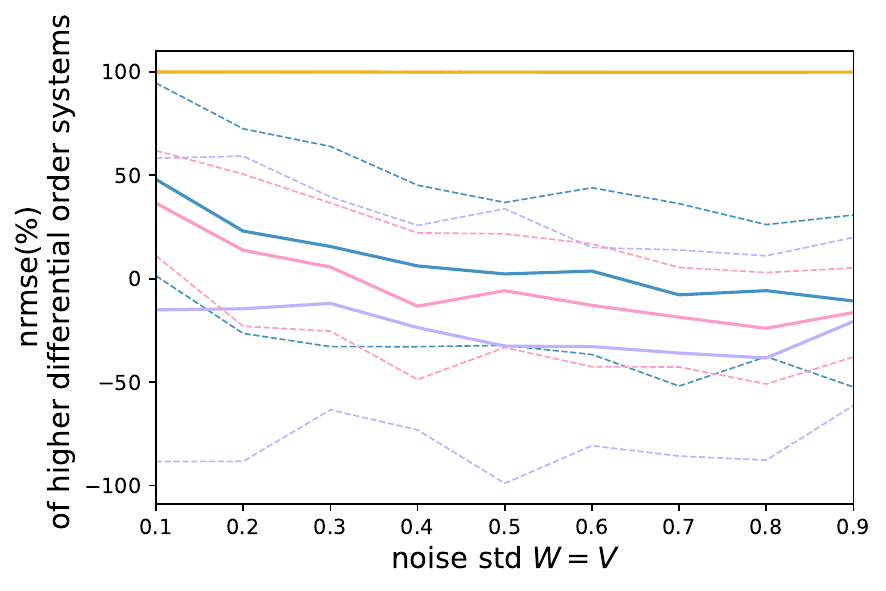}
\includegraphics[width=0.35\textwidth]{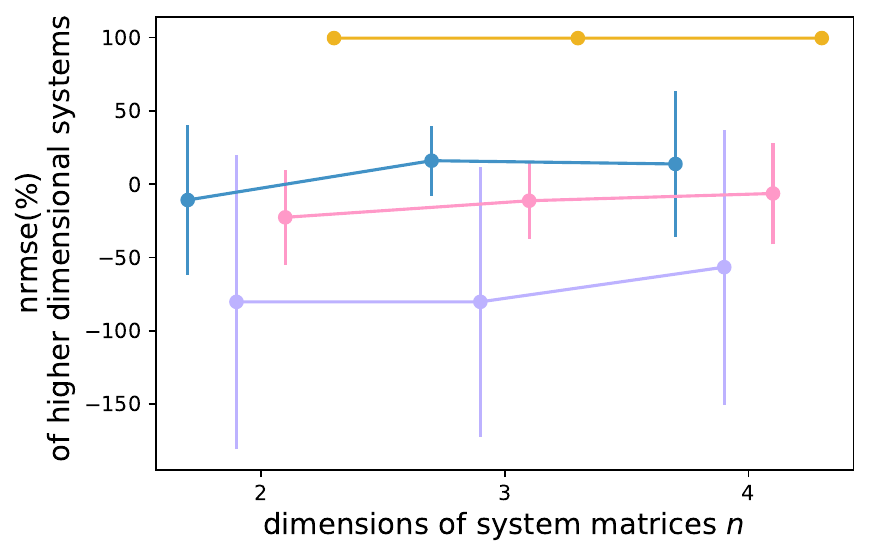}
\caption{The nrmse fitness values \eqref{NRMSE} of our method compared to the leading system identification methods implemented in Matlab\texttrademark{ } System Identification Toolbox\texttrademark.
\textbf{Upper \& middle:} the mean (solid lines) and mean $\pm$ one standard deviations (dashed lines) of nrmse as standard deviation of both process noise and observation noise increasing in lockstep from $0.1$ to $0.9$.
The time series used for simulation are generated from systems in \eqref{equ:LDS} (upper) and higher differential order systems in \eqref{equ:LDS-higher} (middle), with the dimensions $n$ of both systems being $2$. 
\textbf{Lower:} the mean (solid dots) and mean $\pm$ one standard deviations (vertical error bars) of nrmse at different dimensions $n$ of the underlying systems in \eqref{equ:LDS}.
Higher nrmse indicates better simulation performance.
}
\label{fig:CompareSim}
%\caption{higher differential order}
%\caption{higher dimensions of system matrices}
\end{figure}

Figure~\ref{fig:CompareSim} illustrates the results, with methods distinguished by colors: 
blue for ``least squares auto'', purple for ``subspace auto'', pink for ``ssarx'', and yellow for our method.
The upper subplot presents the mean (solid lines) and mean $\pm$ one standard deviations (dashed lines) of nrmse as standard deviation of both process noise and observation noise (``noise std'') increasing in lockstep from $0.1$ to $0.9$. The underlying system is the same as in the upper subplot of Figure~\ref{fig:NCPO100}, except for $W=V=0.1,0.2,\dots,0.9$.
The middle subplot is similar, except the time series are generated by systems of a higher differential order:
\begin{equation}
\begin{split}
\phi_t & = G \phi_{t-1} + \omega_t  \\
Y_{t} & = F_1' \phi_t + F_2'(\phi_t-\phi_{t-1}) + \nu_t,
\end{split}
\label{equ:LDS-higher}
\end{equation}
and the formulation of our method is changed accordingly.
In the lower subplot of Figure~\ref{fig:CompareSim},
we consider the mean (solid dots) and mean $\pm$ one standard deviations (vertical error bars) of nrmse at different dimensions $n=2,3,4$ of the underlying system \eqref{equ:LDS}.

As Figure \ref{fig:CompareSim} suggests, the nrmse values of our method on this example are almost 100\%, while other methods rarely reach 50\%
despite the fact that the dimensions used by the baselines are the true dimensions of the underlying system ($d=n$). 
(We will use ``least squares auto'', which seems to work best within the other methods, in the following experiment on stock-market data.) Additionally, our method shows better stability; the gap between the yellow dashed lines in the upper or middle subplot, which suggests the width of two standard deviations, is relatively small.

\subsection{Experiments with stock-market data}

Our approach to proper learning of LDS could also be used in a ``prediction-error'' method for improper learning of LDS, i.e., forecasting its next observation (output, measurement). As such, it can be applied to any time series. 
To exhibit this, we consider real-world stock-market data first used in \cite{arima_aaai}.
%(The data are believed to be related to the evolution of daily adjusted closing prices of stocks within the Dow Jones Industrial Average, starting on February 16, 1885. Unfortunately, the data are no longer available from the website of the authors of \cite{arima_aaai}.)
In particular, we predict the evolution of the stock price from the 21$^\textrm{st}$ period to the 121$^\textrm{st}$ period, where each prediction is based on the 20 immediately preceding observations ($T=20$). 
%For example, we use the first 20 periods of the stock prices to predict the stock price for the 21$^\textrm{st}$ period. Next, we use the prices from the 2$^\textrm{rd}$ period to the 21$^\textrm{st}$ period to predict the stock price for the 22$^\textrm{st}$ period, and so on.
For our method, 
we use the same formulation \eqref{obj_B} subject to \eqref{NFF_1}-\eqref{NFF_2}, but with the variable $F'$ removed. %We set $c_1=c_2=0.01$ and the moment order $k=1$. 
%The SDP relaxation is generated in \texttt{ncpol2sdpa} 1.12.2, and solved by \texttt{mosek} 9.2.
For comparison,
the combination of least squares algorithms ``least squares auto'' is used again.
Since we are using the stock-market data, the dimension $n$ of the underlying system is unknown.
Hence, the dimensions $d$ of the ``least squares auto'' are iterated from $1$ to $4$, wherein $4$ is the highest setting allowed in the toolbox for $20$-period observations. 

\begin{figure}[t!]
\centering
\includegraphics[width=0.35\textwidth]{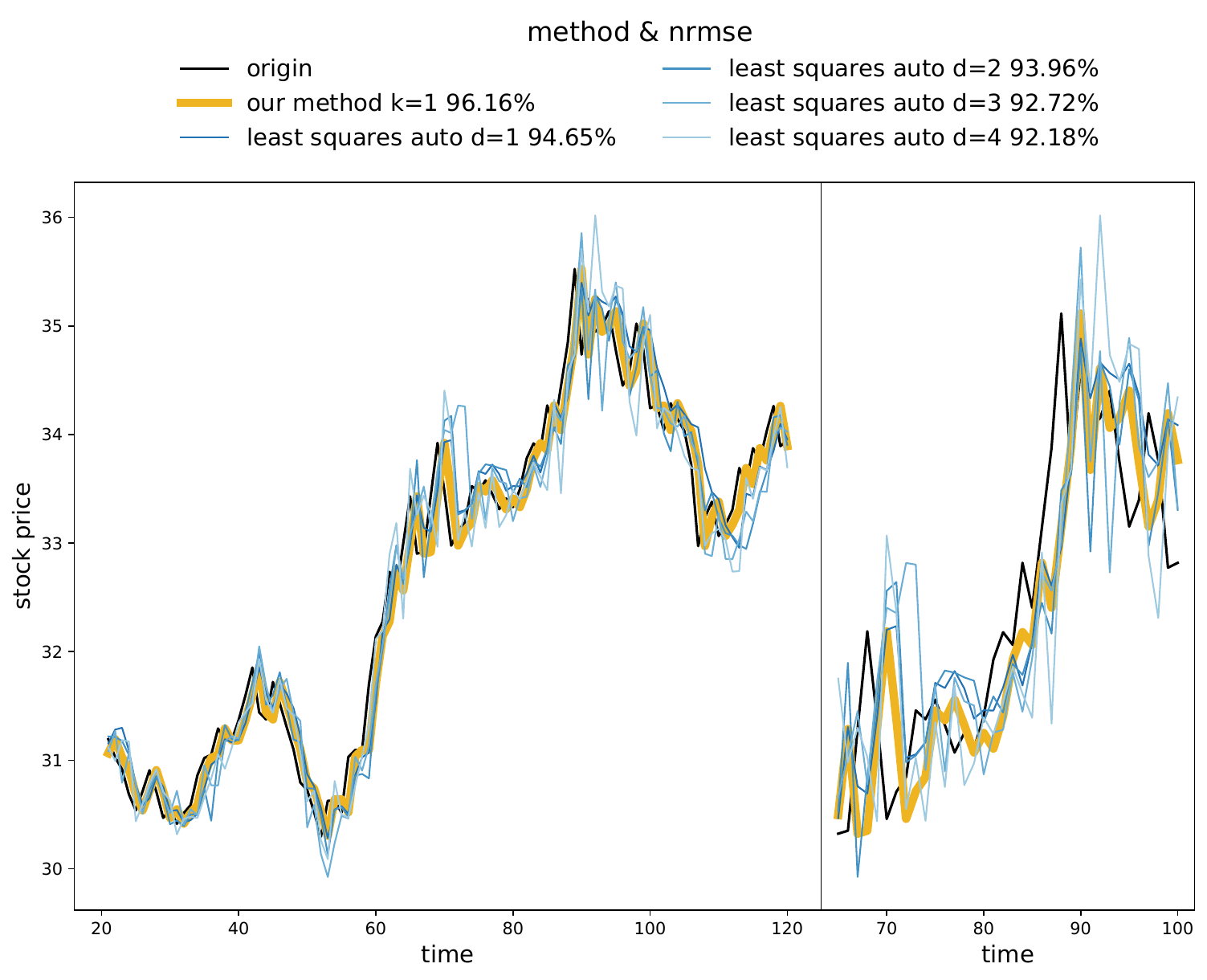}
\caption{\textbf{Left:} 
The time series of stock price (dark) for the 21$^\textrm{st}$-121$^\textrm{st}$ period used in \cite{arima_aaai}, and 
the predicted outputs of our method (yellow) compared against ``least squares auto'' (blue) implemented in Matlab\texttrademark{} System Identification Toolbox\texttrademark.
The dimension $d$ of ``least squares auto'' is iterated from $1$ to the highest number of $4$. 
The percentages in legend are corresponding nrmse values of one-step predictions. \textbf{Right:} a zoom-in for the 66$^\textrm{th}$-101$^\textrm{st}$ period.}
\label{fig:TimePlot}
\end{figure}

Figure~\ref{fig:TimePlot} shows in the left subplot the results obtained by our method (a yellow curve), and the ``least squares auto'' of varying dimensions $d=1,2,3,4$ (four blue curves).
The true stock price ``origin'' is displayed by a dark curve.
%are shown in the left subplot in Figure \ref{fig:TimePlot}. 
The percentages in the legend correspond to nrmse values \eqref{NRMSE}. Both from the nrmse and the shape of these curves, we notice that ``least squares auto'' performs poorly when the stock prices are volatile.
This is highlighted in the right subplot, which zooms in on the 66$^\textrm{th}$-101$^\textrm{st}$ period. 

%which is a zoom-in plot of the 66$^\textrm{th}$ period to the 101$^\textrm{st}$ period. %It seems clear that our method is superior.

\begin{figure}[!t]
\centering
\includegraphics[width=0.4\textwidth]{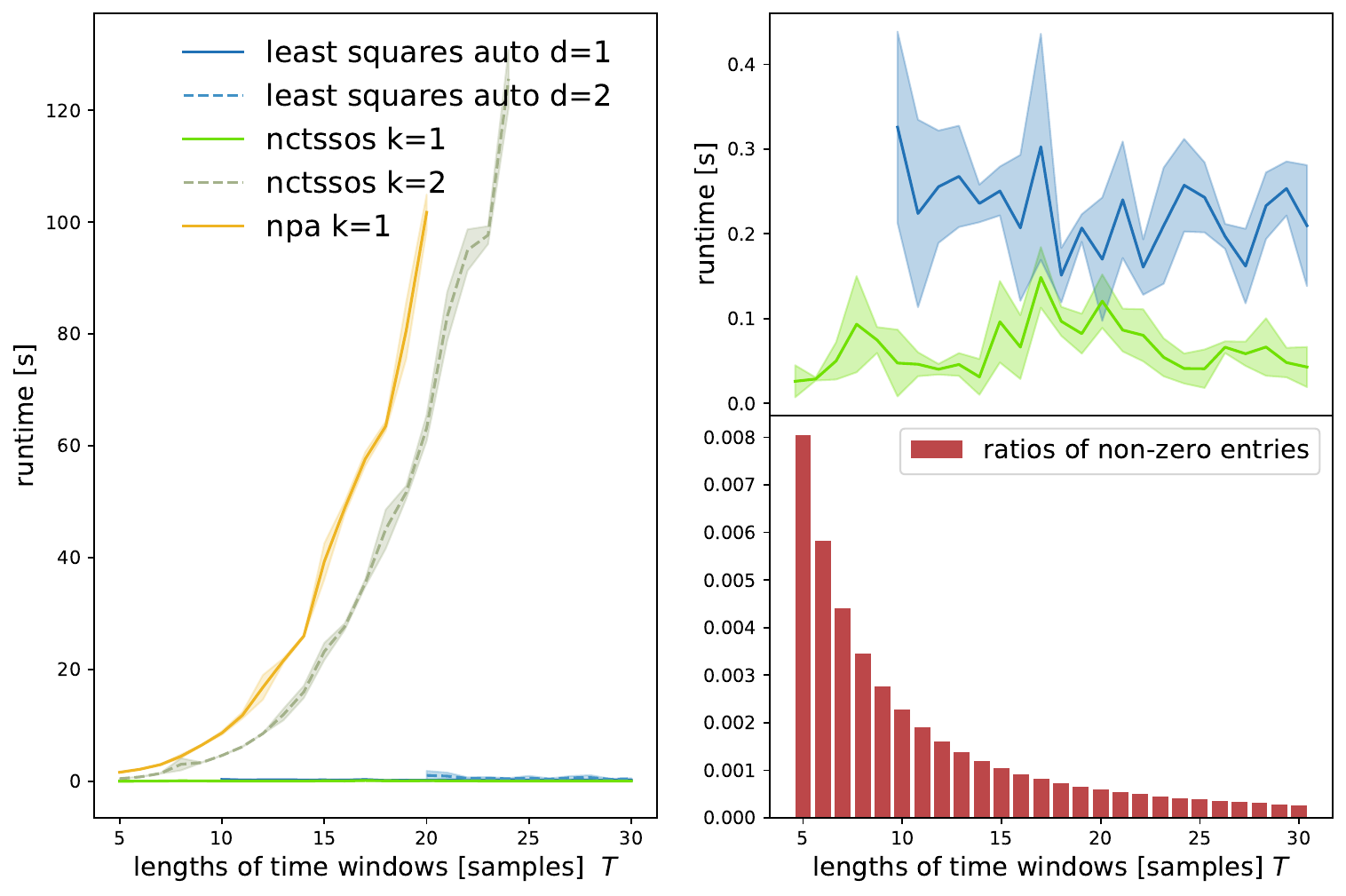}  
\hspace{0in}
\caption{
\textbf{Left:} The (solid or dashed) curves show the mean runtime of the SDP relaxation of the baseline ``least squares auto'' (blue), the TSSOS hierarchy (green) and the NPA hierarchy (yellow), at different moment orders $k$ or dimensions $d$.
The mean $\pm$ one standard deviation of runtime is displayed by shaded error bands.
 \textbf{Upper-right:} The mean and mean $\pm$ one standard deviation of runtime of the SDP relaxation of TSSOS hierarchy at moment order $k=1$ and the ``least squares auto'' with dimension $d=1$.
\textbf{Lower-right:} 
The red bars display the sparsity of NPA hierarchy of the experiment on stock-market data against the length of time window, by ratios of non-zero coefficients out of all coefficients in the SDP relaxations}
\label{fig:runtime}
\end{figure}

%Additionally, we would like to present the run-time of our method as a function of the length of the time window, for univariate observations of the Experiment on Stock-market Data above. Figure \ref{fig:runtime} displays the run-time of our method with varying lengths $T$ of time windows from $5$ to $50$. We can see that the run-time of our method increases fast with the length $T$ of the time window, but that for modest $T$, the absolute run-time may be acceptable in many applications.

%\cite{Pironio2010} of semidefinite programming (SDP) relaxations, as utilized in the proof of Theorem~\ref{T1}, and its sparsity-exploiting variant, known as the term-sparsity exploiting moment/SOS (TSSOS) hierarchy \cite{wang2019tssos,wang2020chordal}.

\subsection{Runtime}

Next, we consider the runtime of two implementations of solvers for \eqref{obj_B} subject to \eqref{NFF_1}-\eqref{NFF_2}. 
The first implementation constructs the SDP relaxation of NPA hierarchy via \texttt{ncpol2sdpa} 1.12.2 with moment order $k=1$. 
The second implementation constructs the non-commutative variant of the TSSOS hierarchy via \texttt{nctssos}, with moment order $k=1,2$. 
%Both relaxations are solved by \texttt{mosek} 9.2.
For comparison purposes, we include the baseline ``least squares auto'' at dimensions $d=1,2$. 
We randomly select a time series from the stock-market data, with the length of time window $T$ chosen from $5,6,\dots,30$, and run these three methods three times for each $T$.
%\textcolor{red}{how to describe each method using different ranges of T.}
 %, i.e., $3\times 3\times 26$ runs in total.

Figure~\ref{fig:runtime} illustrates the runtime of the SDP relaxations and the baseline ``least squares auto'' as a function of the length of the time window. 
These implemented methods are distinguished by colors: blue for ``least squares auto'', green for
the non-commutative variant of the TSSOS hierarchy (``nctssos''), and yellow for the NPA hierarchy (``npa'').
The mean and mean $\pm$ one standard deviation of runtime are displayed by (solid or dashed) curves and shaded error bands.
The upper-right subplot compares the runtime of our method with ``nctssos'' at moment order $k=1$ against ``least squares auto'' with dimension $d=1$.
The red bars in the lower-right subplot display the sparsity of NPA hierarchy of the experiment on stock-market data against the length of time window, by ratios of non-zero coefficients out of all coefficients in the SDP relaxations.

As in most primal-dual interior-point methods \citep{Tuncel2000}, runtime of solving the relaxation to $\epsilon$ error is polynomial in its dimension and logarithmic in $1/\epsilon$, but it should be noted that the dimension of the relaxation grows fast in the length $T$ of the time window and the moment order $k$. 
It is clear that the runtime of solvers for SDP relaxations within the non-commutative variant of the TSSOS hierarchy exhibits a modest growth with the length of time window, much slower than that of the plain-vanilla NPA hierarchy.
%The run-time could be reduced further by combining the exploitation of correlative sparsity \cite{klep2019sparse} and term sparsity \cite{wang2021exploiting}. 
%This should be expected, given that the number of elements in positive-semidefinite constraints in the NPA hierarchy \citep{Pironio2010} is equivalent to the number of entries in the moment matrix $M_k(y)$.
%\textcolor{red}{The sparsity of NPA hierarchy of the experiment is crucial for the significant reduction in runtime when using TSSOS hierarchy.}

\section{Conclusions}

We have presented an alternative approach to the recovery of hidden dynamic underlying a time series, without assumptions on the dimension of the hidden state. 
For the first time in system identification and machine learning, this approach utilizes non-commutative polynomial programming (NCPOP), which has been recently developed within Mathematical Optimization \citep{Pironio2010,wittek2015algorithm,klep2018minimizer,wang2021exploiting}. 
NCPOP can accommodate a variety of other objectives and constraints \cite[e.g. in fairness]{zhou2023fairness}.
This builds upon a long history of work on the method of moments \citep{akhiezer1962some,helton2002positive} and its applications in Machine Learning \citep{LasserreMagron2018}, as well as recent progress \citep{majumdar2019recent} in the scalability of semidefinite programming.

\section*{Acknowledgements}
We would like to thank two anonymous reviewers, Denys Bondar, and Rostislav Horcik for their kind comments, which have helped us improve the draft substantially. 
J.M. acknowledges support of the OP RDE
funded project CZ.02.1.01/0.0/0.0/16\_019/0000765 ``Research Center for Informatics'' and the Czech Science Foundation (23-07947S).
This work has received funding from the European Union’s Horizon Europe research and innovation programme under grant agreement No. 101084642.
This work was also supported by Innovate UK under the Horizon Europe Guarantee; UKRI Reference Number: 101084642 (Codiet).

\bibliographystyle{plain}
\bibliography{ref}

%\clearpage
\onecolumn

%\section*{additional}
%For better clarity of presentation, the experiment in Figure \ref{fig:CompareSim} has been performed again with the same experimental set-up and the results are displayed in Figure \ref{fig:BoxPlot} as four box-plots.
%Each box-plot shows the nrmse values of one method at varying standard deviation of both process noise and observation noise (“noise std”) from $0.1$ to $0.9$ and each box inside indicates the first quartile, mean and the third quartile of 30 nrmse values with the outliers marked by ``o'' symbols. 

%\begin{figure*}[t]
%\centering
%\includegraphics[width=\textwidth]{boxplot_dim2.png}
%\caption{Details of the performance of three other methods and our method used for comparison in Figure~\ref{fig:CompareSim}, in the same experimental set-up.}
%\label{fig:BoxPlot}
%\end{figure*}

\section{Trace Optimization}
\label{app:trace}

Let $\mathbb{S}^{n}$ be the space of $n$-tuple real symmetric matrices of the same order, where $n$ is the number of variables. Suppose $X=\{X_1,\dots, X_n\}\in\mathbb{S}^{n}$, then $X_1,\dots, X_n$ have the same order despite that the specific order is not given.
Let $\mathbb{R}[ X ]_k$ denote the polynomial ring, i.e., polynomials whose coefficients are from $\mathbb{R}$ and whose degrees are less or equal to $k$, or $\infty$ if not specified.
Let $\Sym\mathbb{R}[X]$ denote all symmetric elements in $\mathbb{R}[X]$, such that
\begin{equation}
\Sym\mathbb{R}[X]:=\{p\in\mathbb{R}[X]\;|\;p=p^{\dag}\}.
\label{equ:SymR[X]}
\end{equation}
Let $Q$ be a subset of $\Sym\mathbb{R}[X]$ with elements $q_i, i=1,\dots,m$.
Trace optimization of polynomials in non-commutative variables, with respect to a polynomial $p\in\Sym\mathbb{R}[X]$, has the form in \eqref{trace-NCPO}:
\begin{subequations}
\begin{align}
\tr_{\min}(p)&:=\inf\{\tr p(X)\;|\;X\in\mathbb{S}^{n} \},\label{trace-NCPO-uncon}\\
\tr_{\min}(p,Q)&:=\inf\{\tr p(X)\;|\;X\in\mathbb{S}^{n}; q_i(X)\succcurlyeq 0, i=1,\dots,m\},\label{trace-NCPO-con}
\end{align}
\label{trace-NCPO}
\end{subequations}
where the former formulation is the unconstrained trace minimization, and the latter is the constrained version subject to the positive semidefiniteness of all polynomials in $Q\subseteq\Sym\mathbb{R}[X]$.
Note that if $Q=\varnothing$, then $\tr_{\min}(p,Q)$ becomes an unconstrained problem.

We would need the following definitions to relax these trace minimization problems.

A commutator of two polynomials $p,q\in\mathbb{R}[X]$ is defined as $[p,q]:=pq-qp$.
Since operators can be cyclically interchanged inside a trace, i.e., $\tr(pq)=\tr(qp)$, the trace of a commutator is zero. Conversely, trace zero matrices are (sums of) commutators, using the proof in \cite{albert1957matrices}.
This fact allows us to define an equivalence relation between two polynomials when they have the same trace.
\begin{defin}[\cite{burgdorf2016optimization}, Definition~1.50]{Cyclic equivalence.}
%\paragraph{Cyclic equivalence}
Two polynomials $p,q\in\mathbb{R}[X]$ are cyclic equivalent if $p-q$ is a sum of commutators, denoted by 
\begin{equation}
    p\stackrel{\cyc}{\sim} q.
\end{equation}
\end{defin}
\begin{prop}[\cite{burgdorf2016optimization}, Proposition~1.51]{}
Given two monomials $\omega,\mu$. Then
$\omega\stackrel{\cyc}{\sim}\mu$ if and only if there exist $\nu_1,\nu_2\in\mathcal{W}$, such that $\omega=\nu_1\nu_2$ and $\mu=\nu_2\nu_1$. 
Given two polynomials $p,q\in\mathbb{R}[X]$. Following \eqref{LCoP}, they can be written as $p(X)=\sum_{|\omega|\leq \deg(p)} p_{\omega} \omega$ and $q(X) = \sum_{|\mu|\leq \deg(q)} q_{\mu} \mu$, with $p_{\omega},q_{\mu}\in\mathbb{R}$. Then
\begin{equation}
    p\stackrel{\cyc}{\sim} q \iff \sum_{\omega\in\mathcal{W},\omega\stackrel{\cyc}{\sim}\mu} p_{\omega}=\sum_{\omega\in\mathcal{W},\omega\stackrel{\cyc}{\sim}\mu} q_{\mu},\; \forall \mu\in\mathcal{W}.\label{equ:cyclic equivalence is linear}
\end{equation}
As shown in \eqref{equ:cyclic equivalence is linear}, cyclic equivalence between two polynomials could be formulated as a set of linear constraints.
\label{prop:cyclic equivalence is linear}
\end{prop}

%Recall the definition of polynomial rings from Section~\ref{sec:GNS}: 
%Let $\mathbb{R}[ X ]_k$ denote the polynomial ring, i.e., polynomials whose coefficients are from $\mathbb{R}$ and whose degrees are less or equal to $k$, or $\infty$ if not specified. 
Given a subset of polynomials $Q\subseteq\Sym\mathbb{R}[X]$. With respect to $Q$, recall from Section~\ref{sec:ncpop} that the quadratic module $\mathbf{M}_Q$ is the set of $\sum_if_i^{\dag}f_i+\sum_i\sum_j g_{ij}^{\dag}q_ig_{ij}$, where $f_i$ and $g_{ij}$ are polynomials from the same ring. 
We would use them to define the tracial variants, such that the polynomials within this set have nonnegative trace.
%quadratic modules and the truncated variants.
\begin{defins}[\cite{burgdorf2016optimization}, Definitions~1.22 and~1.56]{Quadratic modules and cyclic quadratic modules.}
%\paragraph{Cyclic quadratic modules}
Given a subset of polynomials $Q\subseteq\Sym\mathbb{R}[X]$.
%With respect to $Q$, recall from Section~\ref{sec:ncpop} that the quadratic module $\mathbf{M}_Q$ is the set of $\sum_if_i^{\dag}f_i+\sum_i\sum_j g_{ij}^{\dag}q_ig_{ij}$, where $f_i$ and $g_{ij}$ are polynomials from the same ring. 
Formally, the quadratic module generated by $Q$ is
\begin{equation}
    \mathbf{M}_Q:=\left\{\sum^{N}_{i=1} f_i^{\dag}f_i+\sum^{m}_{i=1} \sum^{N_i}_{j=1}  g_{ij}^{\dag}q_ig_{ij} \;|\; q_i\in Q;N,N_i\in\mathbb{N};  f_i,g_{ij}\in\mathbb{R}[X]\right\}.
\end{equation}
Note that when $Q=\varnothing$, the quadratic module $\mathbf{M}_Q$, or alternatively $\mathbf{M}_{\varnothing}$ only contains \textit{sums of hermitian squares} polynomials.
As a tracial variant of $\mathbf{M}_Q$, the cyclic quadratic module $\Theta^2_Q$ includes $\mathbf{M}_Q$ and other polynomials that are cyclically equivalent to elements in $\mathbf{M}_Q$:
\begin{equation}
    \Theta^2_Q:=\left\{f\in\Sym\mathbb{R}[X] \;|\; \exists g\in \mathbf{M}_Q: f\stackrel{\cyc}{\sim} g \right\}.
\end{equation}
We define the truncated quadratic module and the truncated cyclic quadratic module of order $2k$ generated by $Q$:
\begin{eqnarray}
\mathbf{M}_{Q,2k}&:=&\mathbf{M}_Q\bigcap\mathbb{R}[X]_{2k},\\
\Theta^2_{Q,2k}&:=&\left\{f\in\Sym\mathbb{R}[X] \;|\; \exists g\in \mathbf{M}_{Q,2k}: f\stackrel{\cyc}{\sim} g \right\}.
\end{eqnarray}
\end{defins}

\begin{lem}[ \cite{helton2002positive}, Lemma~2.1]{}\label{lem:gram matrix}
For a polynomial $p\in\Sym\mathbb{R}[X]_{2k}$, there exist a (nonunique) Gram matrix $G_p$, such that  
\begin{equation}
p=W^{\dag}_{k} G_p W_{k},\; G^{\dag}_p=G_p,
\end{equation}
where $W_k$ is a vector consisting all monomials in $\mathcal{W}_k$.
\end{lem}
\begin{proof}
A symmetric polynomial $p\in\Sym\mathbb{R}[X]_{2k}$ is a weighted sum of symmetrized monomials, such that 
\begin{equation}
p=\sum_{\omega\in \mathcal{W}_{2k}}p_{\omega}( \omega + \omega^{\dag}),\label{equ:sum of symmetrized monomials}
\end{equation}
where $\omega\in \mathcal{W}_{2k}$ could be written as a product of two monomials
\begin{equation}
\omega=u^{\dag} \nu,
\end{equation}
with $u,\nu\in\mathcal{W}_d$. 
%\textcolor{red}{Recall from Section~\ref{app:GNS} that}
Let $\overrightarrow{p}$ denote the vectors of coefficients of $p\in\mathbb{R}[X]$.
Monomials $u,\nu$ could also be regarded as a polynomial, with
\begin{equation}
u=\overrightarrow{u}'W_d,\; \nu=\overrightarrow{\nu}'W_d.
\end{equation}
Since $\omega=\mathcal{W}^{\dag}_d \overrightarrow{u}\overrightarrow{\nu}'W_d$ and $\omega^{\dag}=\mathcal{W}^{\dag}_d \overrightarrow{\nu}\overrightarrow{u}'W_d$,
we show that the symmetrized monomial $\omega+\omega^{\dag}$ could be associated with a symmetric Gram matrix, such that 
\begin{equation}
\omega + \omega^{\dag}=W^{\dag}_d \left(\overrightarrow{u}\overrightarrow{\nu}'+\overrightarrow{\nu}\overrightarrow{u}'\right)W_d.
\label{equ:a symmetrized monomial}
\end{equation}
Substituting \eqref{equ:a symmetrized monomial} for \eqref{equ:sum of symmetrized monomials}, we then obtain a Gram matrix $G_p$ of $p$, with $G_p=\sum_{\omega\in \mathcal{W}_{2k}}p_{\omega}\left(\overrightarrow{u}\overrightarrow{\nu}'+\overrightarrow{\nu}\overrightarrow{u}'\right)$.
Using the fact that the weighted sum of symmetric matrices is also symmetric, we complete the proof.
\end{proof}

\begin{prop}[\cite{burgdorf2016optimization}, Propositions~1.16 and~3.10]{}
\label{prop:unconstrained LMI}
In the unconstrained case ($Q=\varnothing$), the cones $\mathbf{M}_{\varnothing}$ and $\Theta^2_{\varnothing}$ could be formulated as a set of linear matrix inequalities (LMI):
\begin{eqnarray*}
p\in \mathbf{M}_{\varnothing} &\iff& \exists\;G_p\succcurlyeq0: p=W^{\dag} G_p W,\\
p\in \Theta^2_{\varnothing} &\iff& \exists\;\widetilde{G}_p \succcurlyeq0: p\stackrel{\cyc}{\sim} W^{\dag} \widetilde{G}_p W,
\end{eqnarray*}
where $W$ is a vector consisting of all monomials in $\mathcal{W}$. Since Proposition~\ref{prop:cyclic equivalence is linear} has shown that cyclic equivalence is linear, hence the cone $\Theta^2_{\varnothing}$ is also LMI. 
\end{prop}
Especially, $G_p$ is a Gram matrix of $p$ and we call $\widetilde{G}_p$ as a tracial Gram matrix for $p$. 

\begin{prop}[\cite{burgdorf2016optimization}, Proposition~5.7]{}
\label{prop:constrained LMI}
In the constrained case, $\varnothing\neq Q\subset\Sym\mathbb{R}[X]$. Each element $q_{i}\in Q$ could be written as $q_{i}=\sum_{\omega\in\mathcal{W}} q_{i,\omega} \omega$.
Then the cone $\Theta^2_{Q}$ could be formulated as a set of linear matrix inequalities (LMI):
\begin{eqnarray*}
p\in \Theta^2_{Q,2k} &\iff& \exists\; A\succcurlyeq 0,B^{i} \succcurlyeq 0: 
p_{\omega}=\sum_{\mu,\nu\in\mathcal{W}_{k},\mu^{\dag}\nu\stackrel{\cyc}{\sim}\omega}
A_{\mu,\nu} +\sum_i \sum_{\mu,\nu\in\mathcal{W}_{k_i},z\in\mathcal{W}_{\deg q_i},\mu^{\dag}z\nu\stackrel{\cyc}{\sim}\omega} q_{i,z} B^{i}_{\mu,\nu},
\end{eqnarray*}
where $\omega\in\mathcal{W}$, $k_i=\lfloor 
 k-\deg(q_i)/2\rfloor$, and $i\in\{i\;|\;q_i\in Q\}$. The matrix $A$ is of order $|\mathcal{W}_{k}|$ and for each $q_i$, its corresponding matrix $B^{i}$ is of order $|\mathcal{W}_{k_i}|$.
%where $W$ is a vector consisting all monomials in $\mathcal{W}$. Since %Proposition~\ref{prop:cyclic equivalence is linear} has shown that cyclic equivalence is linear, hence the cone $\Theta^2_{\varnothing}$ is also LMI. 
\end{prop}

\begin{defins}[\cite{burgdorf2016optimization}, Definitions~1.23 and~1.59]{The semialgebraic set and the von Neumann semialgebraic set.}
%\paragraph{The semialgebraic set}
Given a subset of polynomials $Q\subseteq\Sym\mathbb{R}[X]$, The semialgebraic set generated by $Q$ is the class of $n$-tuples $X$ of real symmetric matrices of the same order such that $q(X)\succcurlyeq 0$, for all $q\in Q$. Formally,
\begin{equation}
\mathscr{D}_{Q}:=\{X\in\mathbb{S}^n\;|\; \forall q\in Q: q(X)\succcurlyeq 0\}.
\end{equation}
Note that if $Q=\varnothing$, the semialgebraic set is just $\mathbb{S}^n$.
Furthermore, 
%we define the von Neumann semialgebraic set $\mathscr{D}^{\II}_{Q}$ generated by $Q$, following Definition~1.59 of \cite{burgdorf2016optimization}.
let $\mathscr{F}$ be a type $\II$-von Neumann algebra.
Then $\mathscr{D}^{\II}_{Q}$ is the union of the $\mathscr{F}$-semialgebraic set over all type $\II$-von Neumann algebras with separable predual. It has that $\mathscr{D}_{Q}\subseteq\mathscr{D}^{\II}_{Q}$ (cf. Remark~1.61 of \cite{burgdorf2016optimization}).
\end{defins}
The finite von Neumann algebras are introduced because we cannot work on the algebra of all bounded operators on the infinite dimensional Hilbert space, which does not admit a trace.

Following Proposition~1.25 of \cite{burgdorf2016optimization}, it is easy to see that for all $X\in\mathscr{D}_{Q}$ if $f\in \mathbf{M}_Q$, $f(X)\succcurlyeq 0$, thus $\tr f(X)\geq 0$. 
%Using the definition of cyclic quadratic modules, we have that for $X\in\mathscr{D}_{Q}$ if $f\in\Theta^2_{Q}$, then $\tr f(X)\geq 0$.
According to Proposition~1.62 of \cite{burgdorf2016optimization}, we still have that for $X\in\mathscr{D}^{\II}_{Q}$ if $f\in\Theta^2_{Q}$, then $\tr f(X)\geq 0$. 
Using the new notation,
Problem~\eqref{trace-NCPO-con} could be rewritten as finding the trace minimum in the region of $\mathscr{D}_{Q}$. 
Instead of working with $\mathscr{D}_{Q}$ directly, $\tr^{\II}_{\min}(p,Q)$ provides a lower bound for $\tr_{\min}(p,Q)$, and it is also more approachable to solve:
\begin{equation}
    \tr^{\II}_{\min}(p,Q):=\inf\{\tr p(X)\;|\;X\in\mathscr{D}^{\II}_{Q}\}\leq \tr_{\min}(p,Q).\label{trace-NCPO-con-II}
\end{equation}

%\textcolor{red}{Using the GNS construction in Section~\ref{app:GNS}.} 
Consider the finite dimensional case, a symmetric linear functional $L:\mathbb{R}[X]_{2k}\to\mathbb{R}$ (cf. Section~\ref{app:GNS}) is in bijective corresponse to a Hankel matrix $M_k$, indexed by monomials $u,\nu\in\mathcal{W}_{k}$, such that 
\begin{equation}
(M_k)_{u,\nu}:=L(u^{\dag}\nu).
\end{equation}
Given $Q\subseteq\Sym\mathbb{R}[X]$, we define a localizing matrix $M_{k_i}$ for $q_i\in Q$, which is indexed by monomials $u,\nu\in\mathcal{W}_{k-\lceil\deg(q)/2\rceil}$, such that
\begin{equation}
(M_{k_i})_{u,\nu}:=L(u^{\dag}q\nu).
\end{equation}
%Suppose $M_k$ is an extension of another Hankel matrix $H_{\check{L}}$ with $\rk H_{\check{L}}=\rk M_k$, we say that $M_k$ is flat over $H_{\check{L}}$.
%Further, suppose their corresponding linear functionals being $L:\mathbb{R}[X]_{2k}\to\mathbb{R}$ and $\check{L}:\mathbb{R}[X]_{2k-2\delta}\to\mathbb{R}$, we say that $L$ is $\delta$-flat.
\begin{cor}[]{} 
\label{cor:transfer to Hankel}
Given two polynomials $p,q\in\mathcal{R}[X]$ and $p=\sum_{u} p_u u$, $q=\sum_{\nu} q_{\nu} \nu$.
Following from the definition of the Hankel matrix $M_k$, we have that 
\begin{equation}
L(p^{\dag} q)=\sum_{u,\nu}p_u q_{\nu} L(u^{\dag}\nu)=\sum_{u,\nu}p_u q_{\nu} (M_k)_{u,\nu}=\overrightarrow{p}' M_k\overrightarrow{q}.
\end{equation}
%where $\overrightarrow{p},\overrightarrow{q}$ are vectors consisting of all coefficients $p_{u}$ of $p$ and $q_{\nu}$ of $q$.
\end{cor}
\begin{lem}[\cite{burgdorf2016optimization}, Lemma~1.44]{Positivity of $L$ is equivalent to the positive semidefinitness of $M_k$.}  
According to the definition, $L$ is nonnegative on $\mathbf{M}_{Q,2k}$. Hence, using Corollary~\ref{cor:transfer to Hankel}, it is easy to see that
\begin{equation}
M_k\succcurlyeq 0;\; M_{k_i}\succcurlyeq 0,\;\forall q\in Q.
\end{equation}
\end{lem}
In the tracial cases, $L$ also needs to be zero on commutators, hence $M_k$ is invariant under cyclic equivalence, such that for $u,v,w,z\in [X]_{k}$
\begin{equation}
    (M_k)_{u,v}=(M_k)_{w,z},\; \forall\;u^{\dag}v \stackrel{\cyc}{\sim} w^{\dag}z. 
\end{equation}

\begin{defin}[\cite{burgdorf2016optimization}, Remark~1.64]{The dual cones.}\label{defin:dual cones}
We define the dual cone $\left(\Theta^2_{2k}\right)^{\vee}$ to consist of symmetric linear functionals which are nonnegative on $\Theta^2_{2k}$ and zero on commutators.
Likewise, if given $Q\subseteq\Sym\mathbb{R}[X]$, we also give a constrained dual cone. Using the Hankel matrix $M_k$, the dual cone can also be rewritten as a set of LMI.
\begin{eqnarray}
\left(\Theta^2_{2k}\right)^{\vee}:&=&\{L:\mathbb{R}[X]_{2k}\to\mathbb{R}\ \;|\; L \textrm{ linear}, L(p)=L(p^{\dag}),L(p)\geq 0, \forall p\in \Theta^2_{2k}\},\nonumber\\
&\cong&\{M_k\;|\;M_k\succcurlyeq 0, (M_k)_{u,v}=(M_k)_{w,z},\; \forall u^{\dag}v \stackrel{\cyc}{\sim} w^{\dag}z \}.\\
\left(\Theta^2_{Q,2k}\right)^{\vee}:&=&\{L:\mathbb{R}[X]_{2k}\to\mathbb{R}\ \;|\; L \textrm{ linear}, L(p)=L(p^{\dag}),L(p)\geq 0, \forall p\in \Theta^2_{Q,2k}\}.\nonumber\\
&\cong&\{M_k\;|\;M_k\succcurlyeq 0; M_{k_i}\succcurlyeq 0,\;\forall q\in Q; (M_k)_{u,v}=(M_k)_{w,z},\; \forall u^{\dag}v \stackrel{\cyc}{\sim} w^{\dag}z \}.
\end{eqnarray}
\end{defin}

\subsection{Relaxations}

In terms of \eqref{trace-NCPO-uncon}, it is easy to see that $\tr (p-a)(X)\geq 0, \forall X\in\mathbb{S}^n$ if $a\leq\tr_{\min}(p)$. Hence, \eqref{trace-NCPO-uncon} is equivalent to \eqref{trace-NCPO-uncon-sup}. Following the same reasoning, we can rewrite \eqref{trace-NCPO-con-II} as \eqref{trace-NCPO-con-sup}.
\begin{subequations}
\begin{align}
\tr_{\min}(p)&:=\sup\{a\;|\;X\in\mathbb{S}^{n};\tr (p-a)(X)\geq 0 \},\label{trace-NCPO-uncon-sup}\\
\tr^{\II}_{\min}(Q,p)&:=\sup\{a\;|\;X\in\mathscr{D}^{\II}_{Q};\tr (p-a)(X)\geq 0 \}\leq\tr_{\min}(Q,p).\label{trace-NCPO-con-sup}
\end{align}
\label{trace-NCPO-sup}
\end{subequations}

There are two important propositions that would be crucial for the following relaxations.
\begin{prop}[\cite{burgdorf2016optimization}, Proposition~1.62]{}
Let $Q\subseteq\Sym\mathbb{R}[X]$.
If $f\in\Theta^2_{Q}$, then $\tr f(X)\geq 0, \;\forall X\in\mathscr{D}_{Q}$, and further, $\tr f(X)\geq 0,\;\forall X\in\mathscr{D}^{\II}_{Q}$.
\label{prop:1.62}
\end{prop}
\begin{prop}[\cite{burgdorf2016optimization}, Proposition~1.63]{}
Let $Q\bigcup\{f\}\subseteq\Sym\mathbb{R}[X]$.
Suppose $\mathbf{M}_Q$ is Archimedean.
If $\tr f(X)\geq 0,\;\forall X\in\mathscr{D}^{\II}_{Q}$, then $ f+\epsilon\in\Theta^2_{Q}$, for all $\epsilon>0$.
\label{prop:1.63}
\end{prop}
Note that the $\epsilon$ in Proposition~\ref{prop:1.63} could be arbitrarily small but non-zero. In other words, $f$ is very close to $\Theta^2_{Q}$. 

Along with Propositions~\ref{prop:1.62}-\ref{prop:1.63}, the new formulations in \eqref{trace-NCPO-sup} could be relaxed to \eqref{trace-NCPO-sup-Theta}. %Lemma 5.2.
\begin{subequations}
\begin{align}
\tr_{\Theta^2}(p)&:=\sup\{a\;|\;f-a\in\Theta^2_{2k_1}\} \leq \tr_{\min}(p),\label{trace-NCPO-uncon-sup-Theta}\\
\tr^{k_2}_{\Theta^2}(p,Q)&:=\sup\{a\;|\;f-a\in\Theta^2_{Q,2k_2}\} \leq
\tr^{\II}_{\min}(Q,p)\leq\tr_{\min}(Q,p),\label{trace-NCPO-con-sup-Theta}
\end{align}
\label{trace-NCPO-sup-Theta}
\end{subequations}
where the unconstrained problem is relaxed to a single SDP (cf. Corollary~\ref{cor:primal SDP}), with $2k_1=\min_{ g\stackrel{\cyc}{\sim} f}\deg g$ .
In the constrained case, $\tr^{\II}_{\min}(p,Q)$ is relaxed to a hierarchy of SDP (cf. Corollary~\ref{cor:primal SDP}) with $2k_2\geq\min_{ g\stackrel{\cyc}{\sim} f}\deg g$.

To obtain the dual problem of \eqref{trace-NCPO-uncon-sup-Theta}, we conduct the standard Lagrangian duality approach, such that
\begin{eqnarray}
\tr_{\Theta^2}(p)&:=&
\begin{array}{cl}
\sup_{a} & a  \\
\textrm{s.t.} & p-a\in\Theta^2_{2k}
\end{array} \\
&=&\sup_a \inf_{L\in\left(\Theta^2_{2k}\right)^{\vee}} \left\{ a + L(p-a) \right\}
\label{equ:cone-transfer} \\
&\leq &\inf_{L\in\left(\Theta^2_{2k}\right)^{\vee}}\sup_a \left\{ a + L(p-a)\right\} \label{equ:weak duality}\\
&=& \inf_{L\in\left(\Theta^2_{2k}\right)^{\vee}}\sup_a \left\{a + L(p)-a L(1) \right\}\label{equ:linearity of L}\\
&=& \inf_{L\in\left(\Theta^2_{2k}\right)^{\vee}}  \left\{ L(p) +\sup_a \{a(1- L(1))\} \right\}\label{equ:inner problem}\\
&=&\begin{array}{cl}
\inf & L(p)  \\
\textrm{s.t.} & L\in\left(\Theta^2_{2k}\right)^{\vee}\\ & L(1)=1.
\end{array}
=:L_{\Theta^2}(p)
\label{equ:trace-NCPO-dual}
\end{eqnarray}
The equality in \eqref{equ:cone-transfer} uses the Definition~\ref{defin:dual cones} to transfer the cone constraint. The inequality in \eqref{equ:weak duality} utilizes the weak duality, since the variable $a$ is unbounded. The equality in \eqref{equ:linearity of L} is due to the linearity of $L$. Finally, the inner problem in \eqref{equ:inner problem} is bounded only if $L(1)=1$, such that we obatin the dual problem in \eqref{equ:trace-NCPO-dual}.
Note that the same procedures could be applied to \eqref{trace-NCPO-con-sup-Theta} simply by changing $\Theta^2_{2k}$ to $\Theta^2_{Q,2k}$, and the resulted formulation is denoted by $L^{k}_{\Theta^2}(p,Q)$.

\begin{cor}[]{$\tr_{\Theta^2}(p),\tr^{k}_{\Theta^2}(p,Q)$ are SDP.} 
\label{cor:primal SDP}
%For $\tr_{\Theta^2}(p),\tr^{k}_{\Theta^2}(p,Q)$, 
This conclusion follows directly from Propositions~\ref{prop:unconstrained LMI} and~\ref{prop:constrained LMI}, such that the constrains in $\tr_{\Theta^2}(p),\tr^{k}_{\Theta^2}(p,Q)$ are sets of LMI.
\end{cor}

\begin{prop}[]{$L_{\Theta^2}(p),L^k_{\Theta^2}(p,Q)$ are SDP.}
\end{prop}
\begin{proof}
For $L_{\Theta^2}(p)$, we can transfer this problem to the space of Hankel matrix $M_k$.
We first look at the objective function of \eqref{equ:trace-NCPO-dual}.
Since $p\in\Sym\mathbb{R}[X]_{2k}$, the polynomial $p$ can be expressed by a Gram matrix $G_p$, such that $p=\mathbf{W}^{\dag}_{k} G_p \mathbf{W}_{k}$ and $G_p$ is symmetric.
Using the proof of Lemma~\ref{lem:gram matrix}, we have
\begin{eqnarray*}
L(p)&=&L\left(\sum_{\omega\in \mathcal{W}_{2k}}p_{\omega}( \omega + \omega^{\dag})\right)\\
&=&\sum_{\omega\in \mathcal{W}_{2k}}p_{\omega} L( \omega + \omega^{\dag})\\
&=&\sum_{\omega\in \mathcal{W}_{2k}}p_{\omega}\left( L(u^{\dag}\nu) + L(\nu^{\dag} u)\right)\\
&=&\sum_{\omega\in \mathcal{W}_{2k}}p_{\omega}\left( \overrightarrow{u}' M_k\overrightarrow{\nu} + \overrightarrow{\nu}' M_k\overrightarrow{u}\right)\\
&=&\left\langle M_k,\sum_{\omega\in \mathcal{W}_{2k}}p_{\omega}\left(\overrightarrow{u}\overrightarrow{\nu}'+\overrightarrow{\nu}\overrightarrow{u}'\right)\right\rangle \\
&=& \langle M_k,G_p\rangle 
\end{eqnarray*}
Then, using Definition~\ref{defin:dual cones}, the first constraint $L\in\left(\Theta^2_{2k}\right)^{\vee}$ is a set of LMI . 
The last constraint $L(1)=1$ is equivalent to $\langle M_k, E_{1,1}\rangle=1$, where $E_{1,1}$ is a matrix with all entries 0 but $1$ in the $(1,1)$-entry. 
Further, $L^{k}_{\Theta^2}(p,Q)$ only differs from $L_{\Theta^2}(p)$ by the dual cone constraint. By reusing Definition~\ref{defin:dual cones}, we know that $L^{k}_{\Theta^2}(p,Q)$ is also a SDP.
\end{proof}

\subsection{The exploiting correlative sparsity variant}
\label{sec:sparsitytheory}

As shown in the lower-right subplot of Figure~\ref{fig:runtime}, there is much room of sparsity exploiting in the interest of significant reduction in runtime. Here, we give a brief summary of correlative sparsity exploiting trace optimization, and we refer readers to \cite{klep2019sparse} for more details and to \cite{wang2021exploiting} for term sparsity exploiting variant.

The index set of the tuple $X\in\mathbb{S}^n$ is $I=\{1,2,\dots,n\}$. Suppose there is a partition such that $I_1,I_2,\dots, I_r\subseteq I$ and $\bigcup^r_{i=1} I_i=I$. 
Let $X(I_i)\subset X$ denote the set of variables $X_j$ whose indices belong to $I_i$.
Given a sparse objective function $p\in\Sym\mathbb{R}[X]$, we assume that $p$ can be decomposed as
\begin{equation*}
p=p_1 + p_2 + \cdots + p_r,
\end{equation*}
where $p_i\in\Sym\mathbb{R}[X(I_i)]$.
With the assumptions of boundedness of each subset $X(I_i)$ and running intersection property, cf. Assumptions~2,3 and 2.4 in \cite{klep2019sparse}, we obtain the sparse trace minimization formulations:
\begin{defins}[\cite{klep2019sparse}]{Sparse semialgebraic sets.}
Given a subset of polynomials $Q\subseteq\Sym\mathbb{R}[X]$ and the sparse version of $\mathscr{D}^{\spa}_Q$ is defined as
\begin{eqnarray*}
\mathscr{D}^{\spa}_Q:=\mathscr{D}_Q\bigcap \left\{X\in\mathbb{S}^n \;|\;\exists N\in\mathbb{R}_{+}: N-\sum_{j\in I_i} X^2_j\succcurlyeq 0,i=1,2,\dots,r\right\}.
\end{eqnarray*}
\end{defins}
\begin{defins}[\cite{klep2019sparse}]{Sparse cyclic quadratic modules.}
Given a subset of polynomials $Q\subseteq\Sym\mathbb{R}[X]$ and the sparse version of $\mathscr{D}^{\spa}_Q$. We define
\begin{eqnarray*}
\Theta^{\spa}_{Q,2k}&:=&\Theta^{1}_{Q,2k}+\Theta^{2}_{Q,2k}+\cdots+\Theta^{r}_{Q,2k},\\
\Theta^{i}_{Q,2k}&:=&\left\{f\in\Sym\mathbb{R}[X(I_i)] \;|\; \exists g\in \mathbf{M}^{i}_{Q,2k}: f\stackrel{\cyc}{\sim} g \right\}.
\end{eqnarray*}
For unconstrained cases, we can drop the $Q$ in subscriptions, i.e., $Q=\varnothing$.
\end{defins}
Then, by replacing $\Theta^2_{2k_1}$ in \eqref{trace-NCPO-uncon-sup-Theta} with $\Theta^{\spa}_{2k_1}$, we get the sparse formulation $\tr_{\Theta^{\spa}}(p)$. Similarly, we get the sparse formulation for constrained case $\tr^{k_2}_{\Theta^{\spa}}(p,Q)$ by substituting $\Theta^2_{Q,2k_2}$ in \eqref{trace-NCPO-con-sup-Theta} with $\Theta^{\spa}_{Q,2k_2}$.

Alongside with sparse trace minimization, one can also exploit the constant trace property for SDP relaxations, see \cite{mai2021constant}.

%\tr^{k_2}_{\Theta^2}(p,Q)&:=\sup\{a\;|\;f-a\in\Theta^2_{Q,2k_2}\} \leq\tr^{\II}_{\min}(Q,p)\leq\tr_{\min}(Q,p),\label{trace-NCPO-con-sup-Theta}

\section{The Gelfand-Naimark-Segal Construction}
\label{app:GNS}

The Gelfand-Naimark-Segal (GNS) construction essentially produces a *-representation from a positive linear functional of a C*-algebra on a Hilbert space. Under the Archimedean assumption, this method could be applied to non-commutative polynomials, which are not C*-algebras otherwise. 
The GNS-construction view of \cite{klep2018minimizer} starts with the construction of the Hilbert spaces. 
%Let $\mathbb{R}[ X ]_k$ denote the polynomial ring, i.e., polynomials whose coefficients are from $\mathbb{R}$ and whose degrees are less or equal to $k$, or $\infty$ if not specified.
Let $L:\mathbb{R}[ X ]\rightarrow\mathbb{R}$ be a positive linear functional with $L(\mathbf{M}_Q)\geq 0$ and $\mathbf{M}_Q$ denotes the quadratic module of \eqref{NCPO}. There exists a tuple of bounded operators $X=(X_1,\dots,X_n)$ and a normalized vector $\psi$, such that for all $p\in\mathbb{R}[ X ]$
\begin{equation}
    L(p)=\langle\psi ,p(X)\psi \rangle .\label{equ:sesquilinear-form}
\end{equation}
By Cauchy-Schwarz inequality, the sesquilinear form $L$ in \eqref{equ:sesquilinear-form} induces an inner product on the quotient space $\mathbb{R}[ X ]/\mathcal{N}$, where $\mathcal{N}:=\{p\in\mathbb{R}[ X ] | L(p^{\dag}p)=0\}$, as in Section 1.5 of \cite{burgdorf2016optimization}. 

Notice that the moment matrix $M_k$ of moment order $k$ and its extension $M_{\delta},\delta>k$ are positive semidefinite Hankel matrices, where $M_k$ is a short form of $M_k(y)$.
Associated with $L$, we recall from Section~\ref{sec:ncpop} that the element of moment matrix $M_k$ is defined as
\begin{equation}
    M_k(\nu,\omega) = \langle\psi ,\nu^{\dag}\omega\psi \rangle = L(\nu^{\dag}\omega) \quad\forall\; |\nu|,|\omega| \leq k.
\end{equation}
Further, associated with an operator $X_j$, one could define a sub-matrix of $M_{\delta}$, as $M^j_{\delta}$ whose element is
\begin{equation}
    M^j_{\delta}(\nu,\omega)=L(\nu^{\dag}\; X_j \omega) \quad\forall\; |\nu|,|\omega| \leq k. \label{equ:MomentMatrix-sub}
\end{equation}
Let $\overrightarrow{p},\overrightarrow{q}$ be the vectors of coefficients of $p,q\in\mathbb{R}[ X ]$ and $k=\max\{\deg(p),\deg(q)\}$.
We have
\begin{equation}
    L(p^{\dag}q)=\sum_{\nu\leq\deg(p),\omega\leq\deg(q)} p_{\nu}q_{\omega} L(\nu,\omega)
    =\overrightarrow{p}' M_k \overrightarrow{q}
\end{equation}
%where $\overrightarrow{p},\overrightarrow{q}$ be the vectors of coefficients of $p,q\in\mathbb{R}[ X ]$, and $k=\max\{\deg(p),\deg(q)\}$. %Please note that vectors $\overrightarrow{p},\overrightarrow{q},\overrightarrow{x_jp},\overrightarrow{x_jq}$ are of compatible sizes.
Hence, associated to $M_k,M_{\delta}$ are the Hilbert spaces $\mathcal{H}:=\mathbb{R}[X]_{k}/\ker M_k$ and $\mathcal{D}:= \mathbb{R}[X]_{\delta}/\ker M_{\delta}$, 
and $L$ induces an inner product $L(p^{\dag}q)=\langle [p],[q] \rangle_{M_k}=\langle\overrightarrow{p},M_k\overrightarrow{q} \rangle_2$.
Note that $\ker M_k:=\{p\in\mathbb{R}[ X ]_k | \langle \overrightarrow{p},M_k\overrightarrow{p} \rangle_2=0\}$ denotes the null vectors, and $[p]\in\mathcal{H}$ denotes the equivalence class of $p\in\mathbb{R}[ X ]_d$.
If $M_{\delta}$ is a flat extension of $M_{k}$, such that $\textrm{rank}(M_{\delta})=\textrm{rank}(M_k)$, $\mathcal{H}=\mathcal{D}$.
The Hankel property of $M_{\delta}$ implies that 
\begin{equation}
    f(X)[p]=[fp]\in\mathcal{D},
\end{equation}
for every $f,p\in\mathbb{R}[ X ]$ and $\deg f\leq k,\deg f + \deg p\leq \delta$. Consequently, we have for $f\in\mathbb{R}[ X ]_k$
\begin{equation}
    \langle[1], f(X)[1]\rangle_{M_k}=\langle\overrightarrow{1}, M_k \overrightarrow{f}\rangle_2.
\end{equation}
Further, let $M_k=U S U'$ be a singular value decomposition, where $S$ is a positive definite diagonal matrix whose size is $\textrm{rank} M_k$ and $U' U=I$. Hence, the equivalence classes of columns of $U\sqrt{S}^{-1}$ form an orthogonal basis $\mathcal{B(H)}$ of the Hilbert space $\mathcal{H}$.
Since the operator $X_j$ is determined by 
\begin{equation}
 \langle [1],X_j [1] \rangle_{M_k}=\langle\overrightarrow{1}, M_{\delta}\overrightarrow{X_j 1} \rangle_2 = \langle\overrightarrow{1}, M_{\delta}^j\overrightarrow{1} \rangle_2   ,
\end{equation}
the matrix representation of $X_j$ with respect to the orthogonal basis $\mathcal{B(H)}$ is 
\begin{equation}
    \sqrt{S}^{-1}U' M^j_{\delta} U \sqrt{S}^{-1}.
\end{equation}
Thus, one can use the singular value decomposition to recover the solution. 

\section{Numerical illustrations}
\label{sec:numerical}

Let us now present the implementation of the approach using the techniques of non-commutative polynomial optimization  \citep{Pironio2010,burgdorf2016optimization} and to compare the results with traditional system identification methods.
Our implementation is available online \footnote{\url{https://github.com/Quan-Zhou/Proper-Learning-of-LDS}}.
%In our implementation, we make use of a globally convergent Navascués-Pironio-Acín (NPA) hierarchy \citep{Pironio2010} of semidefinite programming (SDP) relaxations, as utilized in the proof of Theorem~\ref{T1}, and its sparsity-exploiting variant, known as the term-sparsity exploiting moment/SOS (TSSOS) hierarchy \citep{wang2019tssos,wang2020chordal}.
%Because the degrees of objective \eqref{obj_B} and constraints in (\ref{NFF_1}--\ref{NFF_2}) are all less or equal to $2$, the moment order $k$ within the respective hierarchy can start from $k=1$. %and increase by $1$ in each iteration. 
%The SDP of a given order in the respective hierarchy is constructed using \texttt{ncpol2sdpa 1.12.2}\footnote{\url{https://github.com/peterwittek/ncpol2sdpa}} of Wittek \citep{wittek2015algorithm} or the \texttt{tssos}\footnote{\url{https://github.com/wangjie212/TSSOS}} of Wang et al. \citep{wang2019tssos,wang2020chordal} 
%and solved by \texttt{mosek 9.2} or \texttt{sdpa} of Yamashita et al. \citep{yamashita2003implementation}.\textcolor{red}{how to cite mosek} 
%Empirically, we plot the run-time as a function of $T$ in Figure~\ref{fig:runtime}.
%Our implementation is available on-line for review purposes and will be open-sourced upon acceptance. 
%It relies on \cite{wittek2015algorithm}, \cite{cafuta2011ncsostools}, and SDPA of \cite{yamashita2003implementation}. As comparison, we also present the runtime of SDP of exploiting term-sparsity (TSSOS), as pioneered by Wang et al. \cite{wang2019tssos,wang2020chordal}.
The values of all parameters in our implementation are summarized in Table~\ref{tab:parameters}.
\begin{table}[t]
\begin{tabular}{|cc|ccccccc|}
\hline
\multicolumn{2}{|c|}{\multirow{2}{*}{\textbf{}}}      & \multicolumn{7}{c|}{\textbf{data}}                                                                                                                                                                                                                                                                                                                                                                                                                                                                                        \\ \cline{3-9} 
\multicolumn{2}{|c|}{}                                & \multicolumn{1}{c|}{$G$}                                                                                  & \multicolumn{1}{c|}{$F'$}                                                                                                                         & \multicolumn{1}{c|}{$V$}                 & \multicolumn{1}{c|}{$W$}                 & \multicolumn{1}{c|}{$\phi'_0$}                                                                  & \multicolumn{1}{c|}{$n$}                  & $T$                   \\ \hline
\multicolumn{1}{|c|}{\multirow{2}{*}{Fig 1}} & upper  & \multicolumn{1}{c|}{\multirow{4}{*}{$\bigl(\begin{smallmatrix} 0.9&0.2\\0.1&0.1\end{smallmatrix}\bigr)$}} & \multicolumn{1}{c|}{\multirow{3}{*}{$\bigl(\begin{smallmatrix} 1&0.8\end{smallmatrix}\bigr)$}}                                                    & \multicolumn{1}{c|}{$0.1,0.2,\dots,0.9$} & \multicolumn{1}{c|}{$0.1,0.2,\dots,0.9$} & \multicolumn{1}{c|}{\multirow{5}{*}{$\bigl(\begin{smallmatrix} 1 & 1 \end{smallmatrix}\bigr)$}} & \multicolumn{1}{c|}{\multirow{4}{*}{$2$}} & \multirow{4}{*}{$20$} \\ \cline{2-2} \cline{5-6}
\multicolumn{1}{|c|}{}                       & lower  & \multicolumn{1}{c|}{}                                                                                     & \multicolumn{1}{c|}{}                                                                                                                             & \multicolumn{2}{c|}{$0.5$}                                                          & \multicolumn{1}{c|}{}                                                                           & \multicolumn{1}{c|}{}                     &                       \\ \cline{1-2} \cline{5-6}
\multicolumn{1}{|c|}{\multirow{3}{*}{Fig 2}} & upper  & \multicolumn{1}{c|}{}                                                                                     & \multicolumn{1}{c|}{}                                                                                                                             & \multicolumn{2}{c|}{$W=V=0.1,0.2,\dots,0.9$}                                        & \multicolumn{1}{c|}{}                                                                           & \multicolumn{1}{c|}{}                     &                       \\ \cline{2-2} \cline{4-6}
\multicolumn{1}{|c|}{}                       & middle & \multicolumn{1}{c|}{}                                                                                     & \multicolumn{1}{c|}{$F_1'=\bigl(\begin{smallmatrix} 1&1\end{smallmatrix}\bigr)$, $F_2'=\bigl(\begin{smallmatrix} 0.2&0.5\end{smallmatrix}\bigr)$} & \multicolumn{2}{c|}{$W=V=0.1,0.2,\dots,0.9$}                                        & \multicolumn{1}{c|}{}                                                                           & \multicolumn{1}{c|}{}                     &                       \\ \cline{2-6} \cline{8-9} 
\multicolumn{1}{|c|}{}                       & lower  & \multicolumn{1}{c|}{random unitary matrices}                                                              & \multicolumn{1}{c|}{$\bigl(\begin{smallmatrix} 1&1\end{smallmatrix}\bigr)$}                                                                       & \multicolumn{2}{c|}{$0.5$}                                                          & \multicolumn{1}{c|}{}                                                                           & \multicolumn{1}{c|}{$2,3,4$}              & $30$                  \\ \hline
\multicolumn{1}{|c|}{Fig 3}                  & left   & \multicolumn{6}{c|}{stock-market data}                                                                                                                                                                                                                                                                                                                                                                                                                                                            & $20$                  \\ \hline
\multicolumn{1}{|c|}{Fig 4}                  & left   & \multicolumn{6}{c|}{stock-market data}                                                                                                                                                                                                                                                                                                                                                                                                                                                            & $5,6,\dots,30$        \\ \hline
\end{tabular}
\begin{tabular}{|cc|cccc|cc|}
%\hline
\multicolumn{2}{|c|}{\multirow{2}{*}{\textbf{}}}      & \multicolumn{4}{c|}{\textbf{our method}}                                                                                                                                             & \multicolumn{2}{c|}{\textbf{baselines}}                                 \\ \cline{3-8} 
\multicolumn{2}{|c|}{}                                & \multicolumn{1}{c|}{$c_1$}                              & \multicolumn{1}{c|}{$c_2$}                              & \multicolumn{1}{c|}{$k$}                  & SDP relaxation       & \multicolumn{1}{c|}{$d$}                  & error                       \\ \hline
\multicolumn{1}{|c|}{\multirow{2}{*}{Fig 1}} & upper  & \multicolumn{1}{c|}{$5\times 10^{-4}$}                  & \multicolumn{1}{c|}{$10^{-4}$}                          & \multicolumn{1}{c|}{\multirow{6}{*}{$1$}} & \multirow{6}{*}{NPA} & \multicolumn{1}{c|}{}                     &                             \\ \cline{2-4} \cline{7-8} 
\multicolumn{1}{|c|}{}                       & lower  & \multicolumn{1}{c|}{$10^{-4},10^{-4}\sqrt{10},\dots,1$} & \multicolumn{1}{c|}{$10^{-4},10^{-4}\sqrt{10},\dots,1$} & \multicolumn{1}{c|}{}                     &                      & \multicolumn{1}{c|}{}                     &                             \\ \cline{1-4} \cline{7-8} 
\multicolumn{1}{|c|}{\multirow{3}{*}{Fig 2}} & upper  & \multicolumn{1}{c|}{\multirow{3}{*}{$5\times 10^{-4}$}} & \multicolumn{1}{c|}{\multirow{3}{*}{$10^{-4}$}}         & \multicolumn{1}{c|}{}                     &                      & \multicolumn{1}{c|}{\multirow{2}{*}{$2$}} & \multirow{3}{*}{simulation} \\ \cline{2-2}
\multicolumn{1}{|c|}{}                       & middle & \multicolumn{1}{c|}{}                                   & \multicolumn{1}{c|}{}                                   & \multicolumn{1}{c|}{}                     &                      & \multicolumn{1}{c|}{}                     &                             \\ \cline{2-2} \cline{7-7}
\multicolumn{1}{|c|}{}                       & lower  & \multicolumn{1}{c|}{}                                   & \multicolumn{1}{c|}{}                                   & \multicolumn{1}{c|}{}                     &                      & \multicolumn{1}{c|}{$=n$}                 &                             \\ \cline{1-4} \cline{7-8} 
\multicolumn{1}{|c|}{Fig 3}                  & left   & \multicolumn{1}{c|}{\multirow{2}{*}{$10^{-2}$}}         & \multicolumn{1}{c|}{\multirow{2}{*}{$10^{-2}$}}         & \multicolumn{1}{c|}{}                     &                      & \multicolumn{1}{c|}{$1,2,3,4$}            & prediction                  \\ \cline{1-2} \cline{5-8} 
\multicolumn{1}{|c|}{Fig 4}                  & left   & \multicolumn{1}{c|}{}                                   & \multicolumn{1}{c|}{}                                   & \multicolumn{1}{c|}{$1,2$}                & NPA,TSSOS            & \multicolumn{1}{c|}{$1,2$}                & simulation                  \\ \hline
\end{tabular}
\caption{An overview of parameters in the implementation. Please refer to Section~\ref{sec:numerical} A for the definitions of parameters.}
\label{tab:parameters}
\end{table}

\subsection{The general settings}

\paragraph{Data generation}
We have generated the time series of observations using linear dynamical systems \eqref{equ:LDS} detailed next, %with the superscript $m=1$,
by specifying the tuple $(G,F,V,W)$ and the initial hidden state $\phi_0$.
We use $n$ to indicate that the time series of observations were generated using $n \times n$ system matrices,
while we use operator-valued variables to estimate these. 
The standard deviations of process noise and observation noise $W,V$ are chosen from $0.1,0.2,\dots,0.9$. Note that $W$ should be of size $n\times n$, while for simplicity, we use $W=0.1$ to refer to $W=0.1\times I_{n}$, where $I_n$ is the $n$-dimensional identity matrix.
%\textcolor{red}{Let $T$ be the length of time window. }
We would take a sequence of measured output from the LDS as the time series $\{Y_t\}_{t=1}^{T}$, with $T$ be the length of the time window.

%Consider a 2-dimensional system ($n=2$), with$G=\bigl(\begin{smallmatrix} 0.99&0\\1&0.2 \end{smallmatrix}\bigr)$, $F'=\bigl(\begin{smallmatrix} 1&0.8\end{smallmatrix}\bigr)$ and the starting point $m_0'=\bigl(\begin{smallmatrix} 1 & 1 \end{smallmatrix}\bigr)$. 

\paragraph{Our formulation and solvers}
%\label{app:numerical}
For our formulation, we use Equations~\eqref{obj_B} subject to (\ref{NFF_1}--\ref{NFF_2}), where we need to specify the values of $c_1$ and $c_2$.
To generate the SDP relaxation of this formulation as in \eqref{NCPO-R}, we need to specify the moment order $k$. Because the degrees of objective \eqref{obj_B} and constraints in (\ref{NFF_1}--\ref{NFF_2}) are all less than or equal to $2$, the moment order $k$ within the respective hierarchy can start from $k=1$.

In our implementation, we use a globally convergent Navascués-Pironio-Acín (NPA) hierarchy \citep{Pironio2010} of SDP relaxations, as utilized in the proof of Theorem~\ref{T1}, and its sparsity-exploiting variant, known as the non-commutative variant of the term-sparsity exploiting moment/SOS (TSSOS) hierarchy \citep{wang2019tssos,wang2020chordal,wang2021exploiting}.
The SDP of a given order within the NPA hierarchy is constructed using \texttt{ncpol2sdpa} 1.12.2\footnote{\url{https://github.com/peterwittek/ncpol2sdpa}} of Wittek \citep{wittek2015algorithm}.
The SDP of a given order within the non-commutative variant of the TSSOS hirarchy is constructed using the \texttt{nctssos}\footnote{\url{https://github.com/wangjie212/NCTSSOS}} of Wang et al. \citep{wang2021exploiting}, with trace minimization implemented.
Both SDP relaxationa are then solved by \texttt{mosek} 9.2 \cite{mosek2020mosek}.
We refer readers to Section~\ref{app:trace} B for a summary of the sparse version of trace minimization.

\paragraph{Baselines}
We compare our method against other leading methods for estimating state-space models, as implemented in MathWorks\texttrademark{ } Matlab\texttrademark{ } System Identification Toolbox\texttrademark. Specifically, we test against a combination of least squares algorithms implemented in routine \texttt{ssest}  (``least squares auto''), subspace methods of \cite{van} implemented in routine \texttt{n4sid} (``subspace auto''), and a subspace identification method of \cite{jansson2003subspace} with an ARX-based algorithm to compute the weighting, again utilized via \texttt{n4sid} (``ssarx''). 

To parameterize the three baselines, we need to specify the dimension $d$ (i.e., order) of the estimated state-space model. We would set $d=n$ directly or alternatively, iterate from $1$ to the highest number allowed in the toolbox when the underlying system is unknown, e.g., in real-world stock-market data.
Then, we need to specify the error to be minimized in the loss function during estimation.
In fairness to the baselines, we use the one-step ahead prediction error when comparing prediction performance and simulation error between measured and simulated outputs when comparing simulation performance.

\paragraph{The performance index}
To measure the goodness of fit between the ground truth $\{Y_t\}_{t=1}^{T}$ (actual measurements) and the noise-free simulated/ predicted outputs $\{F'm_t\}_{t=1}^{T}$, using different system identification methods, we introduce the
\textit{normalized root mean squared error (nrmse)} fitness value: 
\begin{equation}
\mathrm{\nrmse}:=\left(1- \frac{\left\|Y-F'm\right\|^2}{\left\|Y-\mean(Y)\right\|^2}\right)\times 100\%,\label{NRMSE}
\end{equation}
where $Y$ and $F'm$ are the vectors consisting of the sequence $\{Y_t\}_{t=1}^{T}$ and $\{F'm_t\}_{t=1}^{T}$ respectively.
A higher nrmse fitness value indicates better simulation or prediction performance. 

\subsection{Experiments on the example of Hazan et al.}

In our first experiment, we explore the statistical performance of feasible solutions of the SDP relaxation.
%, in comparison with other  system identification methods.
%In order to generate the sequence of observations, we build a linear system (i.e., the ground true) with the underlying dynamic being 2-dimensional. 
In the upper subplot of Figure~\ref{fig:NCPO100},
utilizing the same LDS as in  \cite{hazan2017learning,Jakub},
we consider: $n=2$, $G=\bigl(\begin{smallmatrix} 0.9&0.2\\0.1&0.1\end{smallmatrix}\bigr)$, $F'=\bigl(\begin{smallmatrix} 1&0.8\end{smallmatrix}\bigr)$, the starting point $\phi_0'=\bigl(\begin{smallmatrix} 1 & 1 \end{smallmatrix}\bigr)$, and $T=20$. 
The standard deviations of process noise and observation noise $W,V$ are chosen from $0.1,0.2,\dots,0.9$.
For our method, we set parameters $c_1=5\times 10^{-4},c_2=10^{-4}$, and the moment order $k=1$, then the formulation is relaxed via \texttt{ncpol2sdpa} 1.12.2 and solved by \texttt{mosek} 9.2.
We therefore performed one experiment on each combination of standard deviations of process $W$ and observation noise $V$, i.e., $9\times 9$ runs in total.

In the lower subplot of Figure~\ref{fig:NCPO100}, we have the same settings as in the upper one, except for $W=V=0.5$ and the parameters $c_1,c_2$ being chosen from $10^{-4},10^{-4}\sqrt{10},\dots,1$. We further perform $9\times 9$ experiments on each combination of $c_1$ and $c_2$.

Figure~\ref{fig:NCPO100} illustrates the nrmse values of $81$ experiments of our method in different combinations of standard deviations of process noise $W$ and observation noise $V$ (upper), and another $81$ experiments in different combinations of $c_1$ and $c_2$ (lower). 
%are illustrated in the upper subplot of Figure~\ref{fig:NCPO100}.
%Note that we only use the feasible solutions of the relaxation problems, which may explain some of the ``non-linear'' nature of the plot.
It seems clear the highest nrmse is to be observed for the standard deviation of both process and observation noises close to $0.5$.
While this may seem puzzling at first, notice that higher standard deviations of noise make it possible to approximate the observations by an auto-regressive process with low regression depth \cite[Theorem 2]{Jakub}. 
The observed behavior is therefore in line with previous results \cite[e.g., Figure 3]{Jakub}.

\subsection{Comparisons against the baselines}

Next, we investigate the simulation performance of our method in comparison with other system identification methods, for varying underlying systems. 
In the upper subplot of Figure~\ref{fig:CompareSim},
we still consider the same LDS as in the upper subplot of Figure~\ref{fig:NCPO100}, except for $W=V=0.1,0.2,\dots,0.9$.
For our method, we set the parameters $c_1=5\times 10^{-4},c_2=10^{-4}$, and the moment order $k=1$, then the relaxation is constructed via \texttt{ncpol2sdpa}  and solved by \texttt{mosek} 9.2.
For the newly added baselines, we set $d=n=2$ and minimize the simulation errors.
Our method and the three baselines are implemented $30$ times at each noise standard deviation, i.e., $4\times 30\times 9$ runs in total, with all methods using the same time series.

In the middle subplot of Figure~\ref{fig:CompareSim},
we consider the same setting as in the upper subplot of Figure~\ref{fig:CompareSim}, but another underlying system, with a higher differential order:
\begin{equation*}
\begin{split}
\phi_t & = G \phi_{t-1} + \omega_t  \\
Y_{t} & = F_1' \phi_t + F_2'(\phi_t-\phi_{t-1}) + \nu_t,
\end{split}
%\label{equ:LDS-higher}
\end{equation*}
where $G=\bigl(\begin{smallmatrix} 0.9&0.2\\0.1&0.1 \end{smallmatrix}\bigr)$, $F_1'=\bigl(\begin{smallmatrix} 1&1\end{smallmatrix}\bigr)$, $F_2'=\bigl(\begin{smallmatrix} 0.2&0.5\end{smallmatrix}\bigr)$.
% For our method, we change the constrains in (\ref{NFF_1}--\ref{NFF_2}) accordingly (is this clear enough?)}, and then the new formulation
The relaxation is built via \texttt{ncpol2sdpa} 1.12.2 and solved by \texttt{mosek} 9.2, with $c_1=5\times 10^{-4},c_2=10^{-3}$ and $k=1$.
Our method and the three baselines are implemented $30$ times at each noise standard deviation, i.e., $4\times 30\times 9$ runs in total, with all methods using the same time series.

In the lower subplot of Figure~\ref{fig:CompareSim},
we consider varying dimensions $n=2,3,4$ of the underlying system. Given a dimension $n$, we let $G$ be a random $n$-dimensional unitary matrix, both $F$ and $\phi_0$ be a $n$-dimensional column vector with all entries being $1$, $T=30$, and $W=V=0.5$.
We only take the real part of the measured outputs as the time series.
The settings for our methods and baselines are the same as in the upper subplot of Figure~\ref{fig:CompareSim}, notably the dimensions used by the baselines is the true dimension of the underlying system ($d=n$). 
Our method and the three baselines are run $30$ times for each dimension $n$, i.e., $4\times 30\times 3$ runs in total, with all methods using the same time series.

Figure~\ref{fig:CompareSim} illustrates the simulation performance of our method, compared with three baselines, for different underlying systems.
These methods are distinguished by colors: blue for ``least squares auto'', purple for ``subspace auto'', pink for ``ssarx'' and yellow for our method.
The upper subplot presents the mean (solid lines) and mean $\pm$ one standard deviations (dashed lines) of nrmse as standard deviation of both process noise and observation noise (``noise std'') increasing in lockstep from $0.1$ to $0.9$, with the underlying system being \eqref{equ:LDS}.
The middle subplot reports the same as the upper one. The difference is that the time series used for simulation are generated by higher differential order systems in \eqref{equ:LDS-higher}, and the formulation of our method changed accordingly.
The lower subplot depicts the mean (solid dots) and mean $\pm$ one standard deviations (vertical error bars) of nrmse at different dimensions $n$ of the underlying systems in \eqref{equ:LDS}.
%Throughout, we use the same time series, whose underlying dynamic is 2-dimensional, with standard deviation of both process noise and observation noise (``noise std'') increasing in lockstep from $0.1$ to $0.9$.  
%Under each noise std, we implemented four methods for 30 runs, whose output are 30 nrmse values for each method. We report the corresponding mean and standard deviation of nrmse values in Figure \ref{fig:CompareSim}, where the solid lines and dashed lines indicate mean and mean $\pm$ one standard deviation, respectively. %extreme values respectively with outliers deleted.
%Notice that higher nrmse is better.
%mean in a solid line, dashed line at mean $\pm$ 1 standard deviation.)

As Figure \ref{fig:CompareSim} suggests, the nrmse values of our method on this example are almost 100\%, while other methods rarely reach 50\%. (We will use ``least squares auto'', which seems to work best within the other methods, in the following experiment on stock-market data.) Additionally, our method shows better stability as the gap between the yellow dashed lines in the upper or middle subplot, which suggests the width of 2 standard deviations, is relatively small. 
%The nrmse values of our method fall between 80\% and 100\% while increase swiftly to more than 90\% when noise std gets higher than $0.1$. 

\subsection{Experiments with stock-market data}

Our approach to proper learning of LDS could also be used in a ``prediction-error'' method for improper learning of LDS, i.e., forecasting its next observation (output, measurement). As such, it can be applied to any time series. 
To exhibit this, we consider real-world stock-market data first used in \cite{arima_aaai}.
(The data are believed to be related to the evolution of daily adjusted closing prices of stocks within the Dow Jones Industrial Average, starting on February 16, 1885. Unfortunately, the data are no longer available from the website of the authors of \cite{arima_aaai}.)

In Figure~\ref{fig:TimePlot}, we predict the evolution of the stock price from the 21$^\textrm{st}$ period to the 121$^\textrm{st}$ period, where each prediction is based on its previous 20-period observations ($T=20$). For example, we use the first 20 periods of the stock prices to predict the stock price for the 21$^\textrm{st}$ period. Next, we use the prices from the 2$^\textrm{rd}$ period to the 21$^\textrm{st}$ period to predict the stock price for the 22$^\textrm{st}$ period, and so on.
For our method, we use the same formulation \eqref{obj_B} subject to \eqref{NFF_1}-\eqref{NFF_2}, but with the variable $F'$ removed. We set $c_1=c_2=0.01$ and the moment order $k=1$. 
The SDP relaxation is generated in \texttt{ncpol2sdpa} 1.12.2, and solved by \texttt{mosek} 9.2.
For comparison,
the combination of least squares algorithms ``least squares auto'' is used again.
Since we are using the stock-market data, the dimension $n$ of the underlying system is unknown.
Hence, the dimensions $d$ of the baseline are iterated from $1$ to $4$, as $4$ is the highest number allowed in the toolbox for $20$-period observations. 
%At each iteration, subspace Gauss-Newton least squares search, adaptive subspace Gauss-Newton search, Levenberg-Marquardt least squares search and steepest descent least squares search are tried in sequence and the first descent direction leading to a reduction in estimation cost is used.
Notably, ``least squares auto'' is set to focus on producing a good predictor model, as the one-step ahead prediction error between measured and predicted outputs is minimized during estimation.

%Our method uses the formulation \eqref{obj_B} subject to \eqref{NFF_1}--\eqref{NFF_2} to estimate the one-step ahead state $m_{21}$, which is outputted as the one-step prediction of the stock-price. 
%This formulation is implemented in \texttt{tssos} \citep{wang2020chordal} with moment order $k=1$.
%For comparison, the combination of least squares algorithms \texttt{ssest} (``least squares auto'' in Figure~\ref{fig:CompareSim}) is implemented again, whose dimensions $n$ are iterated from $1$ to $4$, as $4$ is the highest number allowed in the toolbox for $20$-period observations. At each iteration, subspace Gauss-Newton least squares search, adaptive subspace Gauss-Newton search, Levenberg-Marquardt least squares search and steepest descent least squares search are tried in sequence and the first descent direction leading to a reduction in estimation cost is used.
%Note that in this experiment, \texttt{ssest} is set to focus on producing a good predictor model, as the one-step ahead prediction error between measured and predicted outputs is minimized during estimation.
%Since the underlying dynamic is unknown, the order of the estimated state-space model is tried from 1 to 4 ($d=1,\dots,4$), where 4 is the highest order allowed for 20-period input data.

Figure~\ref{fig:TimePlot} shows the prediction results for the 21$^\textrm{st}$-121$^\textrm{st}$ period, using our method (the yellow curve), and the baselines ``least squares auto'' of varying dimensions $d=1,2,3,4$ (four blue curves).
The true stock price ``origin'' is displayed by a dark curve.
The right subplot is a zoom-in plot of the left one for the 66$^\textrm{th}$-101$^\textrm{st}$ period. 
%are shown in the left subplot in Figure \ref{fig:TimePlot}. 
The percentages in the legend correspond to nrmse values \eqref{NRMSE}. Both from the nrmse and the shape of these curves, we also notice that ``least squares auto'' perform poorly when the stock prices are volatile.
This is highlighted in the right zoom-in subplot.
%which is a zoom-in plot of the 66$^\textrm{th}$ period to the 101$^\textrm{st}$ period. %It seems clear that our method is superior.

%Additionally, we would like to present the run-time of our method as a function of the length of the time window, for univariate observations of the Experiment on Stock-market Data above. Figure \ref{fig:runtime} displays the run-time of our method with varying lengths $T$ of time windows from $5$ to $50$. We can see that the run-time of our method increases fast with the length $T$ of the time window, but that for modest $T$, the absolute run-time may be acceptable in many applications.

%\cite{Pironio2010} of semidefinite programming (SDP) relaxations, as utilized in the proof of Theorem~\ref{T1}, and its sparsity-exploiting variant, known as the term-sparsity exploiting moment/SOS (TSSOS) hierarchy \cite{wang2019tssos,wang2020chordal}.

In Figure~\ref{fig:runtime}, we reuse the stock-market data. Our method, using formulation \eqref{obj_B} subject to \eqref{NFF_1}-\eqref{NFF_2}, is implemented twice: one implementation constructs the SDP relaxation of NPA hierarchy via \texttt{ncpol2sdpa} 1.12.2 with $k=1$, the other implementation constructs the non-commutative variant of the TSSOS hierarchy (yellow) via \texttt{nctssos}, with $k=1,2$. Both relaxations are solved by \texttt{mosek} 9.2.
For comparison purposes, we implement the baseline ``least squares auto'' at dimensions $d=1,2$. 
We randomly select a time series from the stock-market data, with the length of time window T chosen from $5,6,\dots,30$, and run these three methods for three times at each $T$.
%\textcolor{red}{how to describe each method using different ranges of T.}
 %, i.e., $3\times 3\times 26$ runs in total.

Figure~\ref{fig:runtime} illustrates the runtime of the SDP relaxations and the baseline ``least squares auto'' as a function of the length of the time window. 
These implemented methods are distinguished by colors: blue for ``least squares auto'', green for our method using 
the non-commutative variant of the 
TSSOS hierarchy (``nctssos''), and yellow for our method using NPA hierarchy (``npa'').
The mean and mean $\pm$ one standard deviation of runtime are displayed by (solid or dashed) curves and shaded error bands.
The upper-right subplot compares the runtime of our method with ``nctssos'' at moment order $k=1$ against ``least squares auto'' with dimension $d=1$.
The red bars in the lower-right subplot display the sparsity of NPA hierarchy of the experiment on stock-market data against the length of time window, as ratios of non-zero coefficients out of all coefficients in the SDP relaxations.

As in most primal-dual interior-point methods \citep{Tuncel2000}, runtime of solving the relaxation to $\epsilon$ error is polynomial in its dimension and logarithmic in $1/\epsilon$, but it should be noted that the dimension of the relaxation grows fast in the length $T$ of the time window and the moment order $k$. 
It is clear that the runtime of solvers for SDP relaxations within the non-commutative variant of the TSSOS hierarchy exhibits a modest growth with the length of time window, much slower than that of the plain-vanilla NPA hierarchy. This should be expected, given that the number of elements in positive-semidefinite constraints in the NPA hierarchy \citep{Pironio2010} is equivalent to the number of entries in the moment matrix $M_k(y)$.
The sparsity of NPA hierarchy of the experiment is crucial for the significant reduction in runtime when using TSSOS hierarchy.
We refer readers to correlative sparsity exploiting variant in \cite{klep2019sparse} and term sparsity exploiting variant in \cite{wang2021exploiting}, or Section~\ref{app:trace} B for a brief summary.

\end{document}